\RequirePackage{ifpdf}
\ifpdf 
\documentclass[pdftex]{sigma}
\else
\documentclass{sigma}
\fi

\numberwithin{equation}{section}

\def \g  {{\mathfrak{g}}}
\def \p  {{\mathfrak{p}}}
\def \sl {\mathfrak{sl}}
\def \Tr { {\mathrm{Tr}} }
\def \Id { {\mathrm{Id}} }
\def \d { {\mathrm{d}} }
\def \w { {\wedge} }

\begin{document}

	
\renewcommand{\PaperNumber}{005}

\FirstPageHeading

\renewcommand{\thefootnote}{$\star$}

\ShortArticleName{Poisson Manifolds, Lie Algebroids, Modular Classes: a Survey}

\ArticleName{Poisson Manifolds, Lie
Algebroids, Modular Classes:\\ a Survey\footnote{This paper is a
contribution to the Proceedings of the 2007 Midwest
Geometry Conference in honor of Thomas~P.\ Branson. The full collection is available at
\href{http://www.emis.de/journals/SIGMA/MGC2007.html}{http://www.emis.de/journals/SIGMA/MGC2007.html}}}

\Author{Yvette KOSMANN-SCHWARZBACH}

\AuthorNameForHeading{Y. Kosmann-Schwarzbach}

\Address{Centre de Math\'ematiques Laurent Schwartz,
\'Ecole Polytechnique, 91128 Palaiseau, France}
\Email{\href{mailto:yks@math.polytechnique.fr}{yks@math.polytechnique.fr}}
\URLaddress{\url{http://www.math.polytechnique.fr/~kosmann/}}

\ArticleDates{Received August 31, 2007, in f\/inal form January
02, 2008; Published online January 16, 2008}

\Abstract{After a brief summary of the main properties of
Poisson manifolds and Lie algebroids in
general,
we survey recent work on the modular classes of Poisson and twisted
Poisson manifolds,
of Lie algebroids with a Poisson or twisted Poisson structure,
and of Poisson--Nijenhuis manifolds. A review of the spinor
approach to the modular class concludes the paper.}


\Classification{17B70; 17B56; 17B66; 53C15; 58A10; 15A66; 17B81; 17-02; 53-02; 58-02}

\rightline{\it Dedicated to the memory of Thomas Branson}

\renewcommand{\thefootnote}{\arabic{footnote}}
\setcounter{footnote}{0}

\pdfbookmark[1]{Introduction}{int}
\section*{Introduction}

The aim of this article is to survey recent work on modular
classes in Poisson and Lie algebroid theory, after brief\/ly
recalling basic def\/initions concerning Poisson manifolds, Lie
algebroids and bialgebroids, and Gerstenhaber algebras.

{\it Poisson
manifolds} are smooth manifolds with a Poisson bracket on their ring of
functions (see~\cite{Vaisman0} and Section~\ref{poisson}).
Since their introduction by
Lichnerowicz in 1977 \cite{Lichnerowicz}, Poisson geometry and
the cohomology of Poisson manifolds have developed into
a wide f\/ield of research.
{\it Lie algebroids} are vector bundles with a Lie bracket on
their space of sections (see \cite{M1,M2,CW} and Section
\ref{alg}).
Fundamental examples of Lie algebroids include both the tangent
bundles of smooth manifolds with the usual Lie bracket of vector
f\/ields,
and f\/inite-dimensional Lie algebras, which
can be considered as Lie algebroids over a point.
Lie algebroids were f\/irst introduced by
Pradines in 1967 \cite{Pradines}
as the inf\/initesimal objects associated with the dif\/ferentiable groupoids
of Ehresmann \cite{Ehresmann} later called {\it Lie groupoids}.
They appear most prominently in Poisson
geometry, as shown by Coste, Dazord and Weinstein \cite{CDW}, in
the theory of group actions and in connection theory \cite{M1,M2},
in foliation theory
\cite{Moerdijk, Fernandes1}, in sigma-models (see, e.g.,
\cite{Bojowald, Kotov, Zabzine, Pestun}), and in many
other situations. Lie--Rinehart algebras, known
under many dif\/ferent names
\cite{Palais,yksM,Huebschmann1},
are algebraic analogues of Lie algebroids.

The def\/inition of the {\it modular class of a
Poisson manifold} f\/irst appeared
in 1985 in the work of Koszul \cite{Koszul}, who did not give it a name,
and, later, in the work of Weinstein \cite{Weinstein}: he demonstrated
the relationship of this notion to the group of {\it modular automorphisms
of von Neumann algebras}, whence the name he proposed, proved many
properties of the {\it modular vector fields} of Poisson manifolds,
and indicated
how to generalize the modular class to Lie algebroids.
The theory was then developed in complex analytic geometry by
Brylinski and Zuckerman~\cite{BZ}, in algebraic geometry by
Polishchuk \cite{Polishchuk},
in the framework of Lie--Rinehart algebras by Huebschmann
\cite{Huebschmann3,Huebschmann2}, and in that of Lie algebroids by
Evens, Lu and Weinstein \cite{ELW} who formally
def\/ined the {\it modular class of a Lie algebroid},
by Ping Xu \cite{Xu1}, and by myself \cite{yks2000}.
The modular class of a Poisson manifold is the obstruction to the
existence of a density invariant under the f\/lows of all Hamiltonian
vector f\/ields. In the case of the Lie algebroid associated with a~foliation with an orientable normal bundle,
the modular class
is the obstruction to the existence of an invariant transverse
measure\footnote{More generally, to each Lie algebroid is associated
a singular foliation whose
leaves are integral manifolds of the image of its
anchor -- a generalized
distribution on the base manifold --,
and the vanishing of the modular class
generalizes the existence of an invariant transverse measure.
Lie algebroids are the setting in which singular foliations can be
studied: the holonomy and the linear holonomy of Lie
algebroids introduced by Fernandes \cite{Fernandes1, Fernandes2}
def\/ined in terms of $A$-paths,
paths whose image under the anchor map is the tangent to their
projection on the base manifold,
generalize the holonomy
and linear holonomy of foliations.}.
In the case of a Lie algebra considered as a Lie algebroid over a
point, the modular class is the trace of
the adjoint representation, the {\it infinitesimal modular
  character}. This is another justif\/ication for the
terminology.

One can consider both the theory of Poisson manifolds
and that of {\it Lie algebras with a~triangular $r$-matrix},
i.e., a
skew-symmetric solution of the classical Yang--Baxter equation, as
particular cases of
the more general notion of {\it Lie algebroid with a Poisson structure}.
To each of these theo\-ries, there corresponds, more generally,
a twisted version.
{\it Twisted Poisson manifolds} were f\/irst considered in some
quantization problems in string theory \cite{Cornalba, Park},
then they were def\/ined as geometric objects
\cite{KStrobl,SW}, and f\/inally viewed as particular cases of Lie
algebroids with twisted Poisson structures \cite{Roytenberg,yksL}.
Just as a Lie algebra with a triangular $r$-matrix
gives rise to a
{\it Lie bialgebra} in the sense of Drinfeld \cite{Drinfeld}, a Lie algebroid with a
Poisson structure gives rise to a {\it Lie bialgebroid} in the sense of
Mackenzie and Xu \cite{MX,yks1995}, while a twisted Poisson structure
gives rise to a {\it quasi-Lie-bialgebroid}\,\footnote{The various
generalizations of Lie bialgebras
were f\/irst considered and studied in \cite{yks1992}, following the
introduction of
the Lie-quasi-bialgebras as the semi-classical limits of quasi-Hopf
algebras by Drinfeld \cite{Drinfeld1989}, while the algebroid
generalization is due to Roytenberg \cite{Roytenberg}.}.
The problem of def\/ining the {\it modular class of a manifold with a
twisted Poisson structure},
and more generally that of a {\it Lie algebroid with a twisted Poisson
structure}, was
solved in \cite{yksL}. (See also a dif\/ferent approach in
the unpublished manuscript \cite{AX}.)

For a Poisson
manifold, the modular class coincides
up to a factor $2$ with the modular class of its cotangent Lie
algebroid (Theorem \ref{theorem1} in Section \ref{pmanifolds}), and
the same is true for a twisted Poisson manifold
(Section \ref{XYZ}), but this simple relation is no longer valid in general for
Poisson or twisted Poisson structures on Lie algebroids, in
particular on Lie algebras.
It was the concept of the {\it modular class
of a morphism of Lie algebroids}, also called a {\it relative modular class},
introduced in \cite{GMM} and \cite{yksW}, that provided a key for
a result valid in this general case (Theorem \ref{untwisted} in
Section~\ref{revisited}, and Section~\ref{XYZ}).
Further results were established in~\cite{yksY} concerning
the case of the regular Poisson and twisted Poisson structures,
whose image remains of constant rank, with
applications to Lie algebra theory (see Section~\ref{regulartriang}).

In Section \ref{section1}, we recall the main def\/initions concerning Poisson
manifolds and Lie algebroids, and we show how the two notions are
connected (Section \ref{1.4}). The modular class of a Poisson manifold
is studied in Section \ref{section2}, that of a Lie algebroid with a Poisson
structure in Section \ref{LAwithpoisson}, with further results in Section
\ref{revisited}.
Section \ref{section4} deals with the case of Lie algebroids in
general, and Theorem \ref{theorem1} relates the modular class of a Poisson
manifold to that of its cotangent Lie algebroid.
In Section \ref{section5}, we introduce the notion of modular
class of a Lie algebroid morphism which enables us to state Theorem
\ref{untwisted}, relating the modular class of a Lie algebroid with
a Poisson structure to that of the Lie algebroid morphism def\/ined by
the Poisson bivector. In Section \ref{regulartriang}, we
summarize results pertaining to the case of regular Poisson
structures, where simplif\/ications occur in the computation of the
modular class (Theorem \ref{theoremregular}),
and we describe the case of quasi-Frobenius and Frobenius Lie
algebras for which the modular class is computed by a simple formula
(Theorem~\ref{frobenius}).
In Section \ref{section7}, we have collected results concerning the
twisted case: the def\/inition of a twisted Poisson structure and
its modular class, the statement of Theorem~\ref{X+Y} in Section~\ref{XYZ}
which gives an explicit
representative of the modular class, the study of the
non-degenerate case with the
linearization of the twisted Poisson condition into the twisted
symplectic condition in Section~\ref{7.4}, and an example of a twisted
triangular $r$-matrix given in Section~\ref{7.5}.
Section \ref{section8} deals with the modular classes of Poisson--Nijenhuis
manifolds and their relation with the bi-Hamiltonian hierarchy of the
Poisson--Nijenhuis structure (Theorem \ref{hierarchy}).
In Section  \ref{section9}, which is the last,
we explain another characterization of the modular class
that has recently appeared, the spinor approach to Dirac structures in
Courant algebroids. Theorem \ref{theoremspin} in Section \ref{spinors}
can be viewed as an alternative
def\/inition of the modular class in the twisted case that conf\/irms
that the original def\/inition in \cite{yksL} is the right one.
In the Appendix, we indicate additional references to recent
publications and work in progress.

Whereas the computation of the
Poisson cohomology ring of a Poisson manifold or the Lie
algebroid cohomology ring of a Lie algebroid is a very dif\/f\/icult
task, there are explicit formulas for computing the modular
class which is the subject of this survey.
The modular class is but the f\/irst in the sequence of higher
characteristic classes of Lie algebroids studied by Kubarski
(see~\cite{Kubarski}), by
Fernandes \cite{Fernandes0,Fernandes1,Fernandes2} and by Crainic
\cite{Crainic}, and
in the framework of generalized gauge theories by Lyakhovich
and Sharapov~\cite{LS}. The generalization of the
properties of the modular class to
the higher order classes is an open problem.

Most, but not all, of the results in this survey are to be found in
several papers written with Camille Laurent-Gengoux, Franco Magri,
Alan Weinstein and Milen Yakimov.
I have tried to give due credit to them and to the other authors whose
work I have consulted. But the literature has become so vast, that I must be
forgiven for not having mentioned all the relevant papers.

\subsection*{Conventions}
All manifolds and bundles are assumed to be smooth, and spaces of
sections, denoted $\Gamma(\cdot)$, are spaces of smooth sections.
Forms (resp., multivectors) on a vector bundle $E$ are
sections of $\w^\bullet E^*$ (resp., sections of $\w^\bullet E$).
We use the notations $\iota$ and $\epsilon$ for the operators of
interior and exterior multiplication on forms and multivectors,
with the convention $\iota_{X\w Y} = \iota_X \circ \iota_Y$.
A derivation of a graded algebra means a graded derivation.
A biderivation of a bracket means a graded derivation in each argument.

\section{Lie algebroids and Poisson geometry}\label{section1}

\subsection{Poisson manifolds}\label{poisson}
Let us recall the def\/inition of a Poisson manifold\footnote{Lie's
{\it function groups}, further studied by Carath\'eodory, were
  precursors of Poisson manifolds.
Poisson manifolds, which occur
already in \cite{FLS},
 were def\/ined by
Lichnerowicz in
  \cite{Lichnerowicz}. Closely related notions also appeared in the
  work of Berezin, Hermann,
Tulczyjew and Kirillov, while Poisson structures on inf\/inite-dimensional
manifolds were introduced under the
name {\it Hamiltonian structures} by several authors working on
integrable systems. For an exposition of Poisson geometry, see
\cite {LM, Vaisman0,M2}.}.
When
$M$ is a smooth manifold
and $\pi$ is a~bivector f\/ield, i.e., a contravariant,
skew-symmetric, $2$-tensor f\/ield, $\pi\in \Gamma (\wedge^{2} TM)$),
def\/ine $\{f,g\}$ as $\pi(\d f,\d g)$, for $f$ and $g \in C^\infty (M)$.
The pair
$(M,\pi)$ is called a {\emph{Poisson manifold}}  if $\{~,~\}$ is an
$\mathbb R$-Lie bracket. When $(M,\pi)$ is a Poisson manifold,
$\pi$ is called a {\it{Poisson
    bivector}} or a~{\it{Poisson structure}} on
$M$ and $\{~,~\}$ is called the {\it{Poisson bracket}}.
Examples are (i) symplectic structures: if $\omega$ is a non-degenerate
closed $2$-form
on $M$, then its inverse is a Poisson bivector, and (ii) the duals of
Lie algebras: if $\g$ is a Lie algebra, then the Poisson bracket on
$\g^*$ is def\/ined by $\{x,y\}(\xi) = \langle\xi,[x,y]\rangle$, where $x$ and $y \in
\g$ are considered as linear forms on $\g^*$, and $\xi \in \g^*$.
This linear Poisson structure on $\g^*$ is called the
{\it Kirillov--Kostant--Souriau Poisson structure}.

We shall call bivector f\/ields simply {\it bivectors}, and, more generally, we
shall call the sections of $\w^\bullet (TM)$ {\it multivectors}.
The space of multivectors is an associative, graded commutative
algebra under the exterior multiplication and it is endowed with a
graded Lie algebra structure (for a shifted grading) called the
{\it Schouten--Nijenhuis bracket} or the {\it Schouten bracket},
making it a~Gerstenhaber algebra (see
Section \ref{gerstenhaber} below). We shall denote the
$(k+\ell-1)$-vector
which is the Schouten--Nijenhuis bracket of a $k$-vector $U$ and an
$\ell$-vector $V$ on $M$ by $[U,V]_M$. The Schouten--Nijenhuis bracket of
multivectors is the extension of the Lie bracket of vector f\/ields as a
biderivation, satisfying $[X,f]_M = X \cdot f$, for $X \in
\Gamma(TM)$, $f \in C^\infty(M)$.
The following fundamental fact is easy to prove.
\begin{proposition}
A bivector $\pi$ is a Poisson structure on $M$ if and only if
$[\pi,\pi]_{M} = 0$.
\end{proposition}

When $\pi$ is a bivector, we def\/ine $\pi^{\sharp} : T^{*}M \to TM$ by
\[
\pi^{\sharp}(\alpha) = \iota_\alpha \pi  ,
\]
where $\alpha \in T^*M$,
and $\iota_\alpha$ denotes the interior product of multivectors by the
$1$-form $\alpha$.

The {\it Hamiltonian vector field} associated to $f \in
C^\infty(M)$ is
\[
H^{\pi}_{f} = \pi^{\sharp} (\d f) .
\]

Associated to any Poisson structure $\pi$ on $M$, there is a
dif\/ferential,  i.e., an operator of degree~$1$ and
square zero, on $\Gamma (\wedge^{\bullet} TM)$,
\[
\d_{\pi} = [\pi , \cdot \, ]_{M}  ,
\]
called the {\it {Lichnerowicz--Poisson differential}} or the
{\it Poisson cohomo\-logy operator}.
The coho\-mo\-logy of the complex  $(\Gamma
(\wedge^{\bullet} TM),
\d_{\pi})$ is called the
{\it Poisson cohomology} of $(M, \pi)$.
In particular, $H^\pi_f = - [\pi,f]_M = - \d_\pi f$.

We also def\/ine the {\it Poisson homology
  operator} on $\Gamma (\wedge^{\bullet} T^{*} M)$ (also called the
{\it Koszul--Brylinski homology operator}) by
\[
\partial_{\pi} = [ \d , \iota_{\pi}]  ,
\]
where the bracket is the graded commutator of operators,
so that $\partial_{\pi} = \d \circ \iota_{\pi} - \iota_{\pi} \circ \d$.
It is an operator of degree $-1$ and square zero.
The homology of the complex  $(\Gamma
(\wedge^{\bullet} T^*M),
\partial_{\pi})$ is called the
{\it Poisson homology} of $(M, \pi)$.

\subsection{Lie algebroids}\label{alg}
We now turn to the def\/inition of Lie algebroids\footnote{For an
exposition of Lie algebroid theory, see \cite{M1, M2} and \cite{CW}.}, and in Section
\ref{1.4}
we shall show how Lie
algebroids appear in Poisson geometry.

A Lie algebroid structure on a real vector bundle $A \to M$ is def\/ined by a
vector bundle map, $a_A : A {\to} TM$,
called the {\it anchor} of $A$,
and an $\mathbb R$-Lie algebra bracket on $\Gamma A$, $[~,~]_{A}$,
satisfying the {\it Leibniz rule},
\[
[X,fY]_{A} = f[X,Y]_{A} + (a_{A} (X) \cdot f) Y  ,
\]
for all $X$, $Y \in \Gamma A$, $f \in C^{\infty}(M)$.
We shall denote such a Lie algebroid by $(A, a_A, [~,~]_A)$ or simply
by $A$, when there is no risk of confusion.
It follows from these axioms that the map $a_{A}$ induces a morphism
of Lie algebras from $\Gamma A$ to $\Gamma (TM)$, which we denote by the
same symbol,  i.e.,
\[
a_{A} ([X,Y]_{A}) = [a_{A} (X) , a_{A} (Y)]_M ,
\]
for all $X$, $Y \in \Gamma A$.
Fundamental examples of Lie algebroids are the following.
\begin{example} Any tangent bundle
$A = TM$ of a manifold, $M$, with $a_{A} = \Id_{TM}$ and the usual Lie bracket of
vector f\/ields, is a Lie algebroid.
\end{example}

\begin{example} Any Lie algebra,
$A = \mathfrak g$, considered as a vector bundle over the singleton $M
= \{\rm{pt}\}$, with trivial anchor $a_A = 0$, is a Lie algebroid.
\end{example}

\begin{example} Lie algebra bundles
(see \cite{M1, M2}) are
locally trivial Lie algebroids with vanishing anchor.
The Lie bracket is def\/ined pointwise and varies smoothly
with the base point. The integration of such bundles of Lie
algebras into bundles of Lie groups
was studied by Douady and Lazard~\cite{DL}.
\end{example}

\begin{proposition}
Associated to a given Lie algebroid,
$(A,a_{A} , [~ ,~ ]_{A})$, there is a
differential, $\d_A$,  on
$\Gamma (\w ^\bullet A^{*})$,
defined by
\begin{gather}
(\d_{A} \alpha) (X_{0},\ldots,X_{k}) = \sum_{i=0}^{k} (-1)^{i} a_{A}
(X_{i}) \cdot  \alpha (X_{0},\ldots, \widehat{X}_{i},\ldots, X_{k}) \nonumber\\
{}  \phantom{(\d_{A} \alpha) (X_{0},\ldots,X_{k}) =}{} +
\sum_{0 \leq i < j \leq k} (-1)^{i+j} \alpha ([X_{i},
X_{j}]_{A}, X_{0},\ldots, \widehat{X}_{i},\ldots,
\widehat{X}_{j},\ldots,X_{k}), \nonumber
\end{gather}
for $\alpha \in \Gamma (\wedge^{k}A^{*})$, $X_{0},\ldots, X_{k} \in
\Gamma A$.
\end{proposition}

The coboundary operator $\d_A$ is called the {\it Lie algebroid
  differential} of $A$. The cohomology of the complex
$(\Gamma (\wedge^\bullet A^{*}), \d_A)$,
called the {\it {Lie algebroid cohomology}}
of $A$,
unif\/ies de Rham and Chevalley--Eilenberg cohomologies \cite{Palais}
(see also \cite{Huebschmann1}): for a tangent bundle,
$A = TM$, with $a_{A} = \Id_{TM}$, $\d_{A} = \d_M$ is the de Rham
dif\/ferential, while for a Lie algebra,
$A = \mathfrak g$, considered as a Lie algebroid over the singleton $M
= \{\rm{pt}\}$, $\d_{A} = \d_{\g}$ is the Chevalley--Eilenberg dif\/ferential.

We call sections of $\w ^\bullet A$ (resp., $\w^\bullet (A^*)$)
{\it multivectors} (resp., {\it forms}) on $A$.

\subsection[Gerstenhaber algebras and Batalin-Vilkovisky algebras]{Gerstenhaber algebras and Batalin--Vilkovisky algebras}\label{gerstenhaber}

Gerstenhaber showed that there is a graded Lie bracket, since called a
{\it Gerstenhaber bracket},
on the Hochschild cohomology of an
associative algebra \cite{Gerstenhaber}.
We recall the def\/inition of a Gerstenhaber algebra.

A {\it Gerstenhaber algebra} is an associative,
graded commutative algebra over a f\/ield $K$ of
cha\-racte\-ristic zero, $(\mathcal A = \oplus_{k\in \mathbb Z}
\mathcal A ^k, \w)$,
with a graded $K$-Lie algebra
structure, $[~,~]$, on $(\oplus_{k\in \mathbb Z}
\mathcal A ^{k+1}, \w)$
such that,  for $u \in \mathcal A ^k$,
$[u, \cdot \, ]$ is a derivation
of degree $k-1$ of $(\mathcal A = \oplus_{k\in \mathbb Z}
\mathcal A ^k, \w)$,
\[
[u, v \w w] = [u,v] \w w + (-1)^{(k-1)\ell} v
\w [u,w]  ,
\]
for $u \in \mathcal A ^k$, $v \in \mathcal A ^\ell$, $w \in \mathcal A$.
The bracket of a
Gerstenhaber algebra is also called an {\it odd Poisson bracket}
\cite{Vaintrob,Voronov} or a {\it Schouten bracket} \cite{yksM}.

A linear map,
$\partial : \mathcal A \to \mathcal A$, of degree $-1$,
is called a {\it generator of the Gerstenhaber bracket} if
\begin{equation}\label{generator}
[u,v]_{\mathcal A} = (-1)^{k} (\partial(u\wedge v) - \partial u \wedge v -
(-1)^{k} u \wedge \partial v) ,
\end{equation}
for all $u\in {\mathcal A}^k$, $v \in {\mathcal A}$.
If $\partial^{2} = 0$, then $({\mathcal A},\partial)$
is called a
{\it Batalin--Vilkovisky algebra} (in short, BV-{\it algebra}).

\medskip

For a Lie algebroid, $A$, there is a Gerstenhaber bracket on $\Gamma
(\wedge^{\bullet} A)$,
which generalizes the Schouten--Nijenhuis bracket of multivectors on a
manifold:

\begin{proposition}
When $A$ is a Lie algebroid,
$(\Gamma (\wedge^{\bullet} A) , \w, [~,~]_{A})$,
where $[~,~]_{A}$ is the extension of the Lie
bracket on $\Gamma A$ as a biderivation satisfying $[X,f]_A = \langle
X,\d_A f \rangle$, for $X \in \Gamma A$ and $f \in C^\infty (M)$,
is a Gerstenhaber algebra.
\end{proposition}

When $A= TM$, for $M$ a manifold, $\Gamma (\w ^\bullet A)$
is the exterior algebra of
multivectors on $M$ and the
Gerstenhaber bracket on this algebra is the Schouten--Nijenhuis
bracket, which we have introduced in Section \ref{poisson} and denoted
there by  $[~,~]_M$. When $A= \g$, a Lie algebra,
we recover the {\it algebraic Schouten bracket}
on $\w ^\bullet \g$, denoted $[~,~]_\g$,
which can be def\/ined as the bracket induced from
the Schouten--Nijenhuis bracket of left-invariant multivectors on a
Lie group with Lie algebra $\g$.

In both cases, the Gerstenhaber bracket is a Batalin--Vilkovisky
bracket: for a manifold,
$M$, transposing the de Rham dif\/ferential by means of the isomorphism
from forms to multivectors def\/ined by a top degree density yields a
generator of square zero of the Schouten--Nijenhuis bracket, while the
Lie algebra homology operator yields a generator of square zero of the
algebraic Schouten bracket of a Lie algebra.

The problem of identifying which Gerstenhaber algebras are
Batalin--Vilkovisky algebras is the subject of much
recent work. See
\cite{Menichi} and references therein.

\subsection{The cotangent bundle of a Poisson manifold is a Lie
  algebroid}\label{1.4}
Let $M$ be a manifold equipped with a bivector $\pi\in \Gamma
(\wedge^{2} TM)$.
Def\/ine
\begin{equation}\label{FMMbracket}
[ \alpha ,\beta ]_{\pi} = \mathcal{L}_{\pi^{\sharp}\alpha}
\beta - \mathcal{L}_{\pi^{\sharp}\beta} \alpha - \d (\pi (\alpha
,\beta)) ,
\end{equation}
for $\alpha$, $\beta \in \Gamma (T^{*}M)$, where $\mathcal L$ denotes
the Lie derivation.
It is easy to show that this skew-symmetric bracket\footnote{Formula
  (\ref{FMMbracket}) appeared as the bracket of $1$-forms on a
symplectic manifold in the f\/irst edition of the
book on mechanics by Abraham and
Marsden \cite{AM} in 1967, before anyone spoke about Poisson manifolds,
and again in \cite{AM2}; formula \eqref{dbracket} below is also
stated there. Formula
  (\ref{FMMbracket}) was re-discovered as the bracket of $1$-forms
on a Poisson manifold by several authors, Fuchssteiner
\cite{Fuchssteiner}, Gelfand and Dorfman \cite{GD}, Magri and
Morosi \cite{MM} and Karas\"ev \cite{Karasev}. In 1987, Weinstein made
the connection between this formula and the theory of Lie groupoids
and Lie algebroids (see \cite{CDW}):
the manifold of units $\Gamma_0$ of a
symplectic groupoid $\Gamma$ is a Poisson manifold, and the
Lie algebroid of such a Lie groupoid is canonically isomorphic to the
cotangent bundle of $\Gamma_0$, so
there is a Lie algebroid bracket on the space of sections of the
cotangent bundle of $\Gamma_0$. A calculation proves that it is given
by \eqref{FMMbracket}, the formula of Abraham and Marsden.}
satisf\/ies
the Jacobi identity if and only if $[\pi,\pi]_{M} = 0$.
Bracket (\ref{FMMbracket}) satisf\/ies
\begin{gather}\label{dbracket}
[\d f, \d g]_\pi = \d \{ f, g \}  ,
\end{gather}
for all $f$, $g \in C^\infty(M)$,
and it satisf\/ies the Leibniz rule
\begin{equation}
[\alpha, f \beta]_\pi = f [\alpha, \beta]_\pi + (\pi^\sharp \alpha
\cdot f) \beta ,
\end{equation}
for all $\alpha,\beta \in \Gamma (T^{*}M)$ and $f \in C^\infty(M)$.
It was remarked by Huebschmann \cite{Huebschmann1}
that these properties determine
the above bracket.

The extension to forms of all degrees on a Poisson manifold of bracket
(\ref{FMMbracket})
is denoted by the same symbol. It is
the {\it {Koszul bracket}} \cite{Koszul}, which satisf\/ies the relation
\begin{equation}\label{koszul}
[\alpha, \beta]_\pi =
 (-1)^{k} (\partial_\pi(\alpha \wedge \beta) - \partial_\pi \alpha \wedge \beta -
(-1)^{k} \alpha \wedge \partial_\pi \beta) ,
\end{equation}
for all $\alpha\in \Gamma(\w ^k T^*M)$, $\beta \in \Gamma(\w ^\bullet T^*M)$,
where $\partial_{\pi}$ is the Poisson homology operator
$[ \d , \iota_{\pi}]$. In particular,
$[\alpha, f]_\pi = (\pi^\sharp \alpha) \cdot f$, for $\alpha \in
\Gamma(T^*M)$ and $f \in C^\infty(M)$.
The following proposition follows from \cite{Koszul,BV,yksM}.

\begin{proposition}
If $(M,\pi)$ is a Poisson manifold,
$(T^{*} M , \pi^{\sharp} , [~,~]_{\pi})$ is a Lie algebroid.
The associated Lie algebroid differential
is the Lichnerowicz--Poisson differential on the multivectors on $M$,
$\d_{\pi} = [\pi , \cdot \, ]_{M}$.
The associated Gerstenhaber algebra is the algebra of forms on $M$
with the Koszul bracket, which is a Batalin--Vilkovisky
algebra with generator
the Poisson homology operator $
\partial_{\pi} = [ \d , \iota_{\pi}]$.
\end{proposition}

In fact, more can be said. There is a compatibility condition between
the Lie algebroid structures of $TM$ and $T^*M$ making the pair
$(TM,T^*M)$ a Lie bialgebroid, a notion which we now def\/ine.

\subsection{Lie bialgebroids}
A {\it {Lie bialgebroid}} \cite{MX,yks1995} is a
pair of Lie algebroids in duality $(A,A^*)$ such that
$\d_A$ is a derivation of $[~,~]_{A^*}$.
This condition is satisf\/ied if and only if
$\d_{A^*}$ is a derivation of $[~,~]_{A}$.
Therefore the notion of Lie bialgebroid is self-dual:
$(A,A^*)$ is a Lie bialgebroid if and only if $(A^*,A)$ is a Lie
bialgebroid.

When $\pi$ is a Poisson bivector on $M$ and $T^*M$ is equipped with
the Lie algebroid structure def\/ined above, the pair $(TM,T^*M)$ is
a Lie bialgebroid:
in fact $\d_{T^*M} = \d_{\pi} = [\pi , \cdot \, ]_{M}$ is
clearly a derivation of the Gerstenhaber bracket, $[~,~]_M$.

When $M$ is a point, the notion
of Lie bialgebroid reduces to that of
a {\it Lie bialgebra} \cite{Drinfeld} (see also \cite{LR,yks1992}).

\subsection{Lie algebroids with a
Poisson structure}\label{triangular}
Lie algebroids with a Poisson structure generalize Poisson
manifolds.

Let $(A,a_A,[\ ,\ ]_{A})$ be a Lie algebroid,
and assume that $\pi \in \Gamma (\wedge^{2} A)$ satisf\/ies
$[\pi,\pi]_A =0$. Then $(A,\pi)$ is called a
{\it Lie algebroid with a Poisson structure}.

Let us def\/ine
\begin{equation} \label{FMM}
[ \alpha ,\beta ]_{\pi} = \mathcal{L}^A_{\pi^{\sharp}\alpha}
\beta - \mathcal{L}^A_{\pi^{\sharp}\beta} \alpha - \d_A (\pi (\alpha
,\beta))  ,
\end{equation}
for $ \alpha,  \beta  \in  \Gamma (A^*) $,
where $\mathcal L^A$ denotes
the Lie derivation on sections of $\w ^\bullet (A^*)$ def\/ined by
${\mathcal L}^A_X = [\d_A, \iota_X]$, for $X \in \Gamma A$, and
$\pi^\sharp:A^* \to A$ is def\/ined by $\pi^\sharp \alpha = \iota_\alpha
\pi$, and set
$a_{A^*}= a_A \circ \pi^\sharp$. Then
$(A^*, a_{A^*}, [~,~]_\pi)$ is a Lie algebroid. The associated
dif\/ferential is $\d _\pi = [\pi, \cdot \, ]_A$, and the Gestenhaber
bracket
on the sections of $\wedge^\bullet(A^*)$, also denoted by
$[~,~]_\pi$,
satisf\/ies
relation \eqref{koszul},
where $\partial_\pi = [\d_A, \iota_\pi]$.
It is therefore clear that, when $A$ is a Lie algebroid with a Poisson
structure,
the pair $(A,A^*)$ is a Lie bialgebroid.

The Poisson cohomology
$H^\bullet(A,\d_\pi)$ of a Lie algebroid with a Poisson
structure $(A,\pi)$
is the cohomology of the
complex $(\Gamma(\wedge^{\bullet} A), \d_{\pi})$ which is the
cohomology of the Lie
algebroid $A^*$ with anchor $a_A  \circ \pi^\sharp$ and
Lie bracket $[~,~]_\pi$, and which generalizes the
Poisson cohomology of Lichnerowicz for Poisson manifolds
\cite{Lichnerowicz}.
The Poisson homology
$H_\bullet(A,\partial_\pi)$ of $(A,\pi)$ is that of the complex
$(\Gamma(\wedge^{\bullet}A^*),\partial_{\pi})$ which has been
studied by Huebschmann
\cite{Huebschmann2}, and which generalizes the
Poisson homology of Poisson manifolds \cite{Koszul,Brylinski}.

If $A = \mathfrak g$, a Lie algebra, and $\pi = r\in \w ^2
  \g$, the Poisson condition is
\begin{gather} \label{cybe}
[r,r]_\g=0  .
\end{gather}
This equation is the
{\it{classical Yang--Baxter
equation}} (CYBE)
which can also be written in tensor notation (\cite{Drinfeld} and
see \cite{CP, yksLNP}),
\[
[r_{12} ,  r_{13}] + [ r_{12} , r_{23}] +  [ r_{13} , r_{23}] = 0 .
\]
If $ r\in \w ^2 \g$ satisf\/ies the classical
Yang--Baxter equation, it is
called a {\it{triangular $r$-matrix}}.
As a~particular case of the study of Lie algebroids
with a Poisson structure, we recover the fact that a
triangular $r$-matrix on $\g$ def\/ines a Lie algebra structure on $\g ^*$, by
\begin{equation} \label{rmatrix}
[\alpha, \beta]_r = {\rm{ad}}^*_{r^\sharp \alpha}\beta -
{\rm{ad}}^*_{r^\sharp \beta}\alpha  ,
\end{equation}
in fact a Lie bialgebra structure
on $(\g,\g^*)$, called a {\it{triangular Lie bialgebra}}.

In analogy with the case of
Lie bialgebras def\/ined by a skew-symmetric solution of the classical Yang--Baxter
equation,
a Lie bialgebroid def\/ined by a Poisson bivector
is also called a {\it {triangular Lie bialgebroid}} \cite{LX}.

\medskip
\subsection{The big bracket}
This section is a very brief summary of parts of several articles
\cite{Vaintrob,Roytenberg,Terashima,yksderived}.
The term ``big bracket'' was coined in \cite{yks1992}.

Let $A$ be a vector bundle over $M$.
The {\it big bracket} is the canonical Poisson bracket on the cotangent
bundle of the supermanifold $\Pi A$, where $\Pi$ denotes the change
of parity
\cite{Vaintrob,Voronov}, generalizing the bracket of
\cite {KS,LR,yks1992}.
We shall denote it by $\{~,~\}$.

Let  $\mu$ denote a Lie
algebroid structure of $A$ viewed as a function on $T^*(\Pi A)$. It is
a cubic element of
bidegree $(1,2)$, which is {\it homological}, i.e., $\{\mu,\mu\}
=0$. Then
\[
[X,Y]_A = \{\{X, \mu\}, Y\}  ,
\]
for all $X$, $Y \in \Gamma (\w ^\bullet A)$.
Thus the
Lie algebroid bracket $[~,~]_A$ is a {\it{derived bracket}} of the big
bracket in the sense of \cite{yks1996,yksderived}.

A Lie bialgebroid is def\/ined by a homological function with terms of
bidegrees $(1,2)$ and~$(2,1)$.

A {\it {Lie-quasi-bialgebroid}}
(resp., {\it {quasi-Lie-bialgebroid}}) is def\/ined by a homological cubic
element with no term of bidgree $(3,0)$ (resp., no term of bidegree $(0,3)$).
Whereas for a Poisson manifold~$M$, the pair $(TM,T^*M)$ is a
bialgebroid, it is a quasi-Lie-bialgebroid when $M$ is a~twisted
Poisson manifold (see Section~\ref{7.1}).

A {\it {proto-bialgebroid}}
is def\/ined by a homological cubic element in the algebra of functions on
$T^*(\Pi A)$  \cite{Roytenberg} (for the case of Lie algebras, see~\cite{yks1992}).

The bracket \eqref{FMM} on $\Gamma(A^*)$ def\/ined by a Poisson bivector $\pi$ is
$\{\pi, \mu\}$, and the Gerstenhaber bracket $[~,~]_\pi$ on $\Gamma (\w
A^*)$ satisf\/ies
\[
[\alpha, \beta]_\pi = \{\{\alpha, \gamma \}, \beta\}  ,
\]
for all $\alpha, \beta \in \Gamma (\w^\bullet A^*)$, where $\gamma =
\{\pi, \mu\}$. Thus the
Gerstenhaber bracket $[~,~]_\pi$ also is a {\it{derived bracket}} of the big
bracket.

\section{The modular class of a Poisson manifold}\label{section2}

The modular class of a Poisson manifold is a class in the f\/irst
Poisson cohomology space of the manifold, i.e., the equivalence class
modulo Hamiltonian vector f\/ields of a vector f\/ield which is an
inf\/initesimal Poisson automorphism, called a modular vector f\/ield.
On an orientable manifold, there exists a volume form invariant under
all Hamiltonian vector f\/ields if and only if there exists a modular
vector f\/ield which vanishes.

\subsection[Modular vector fields and modular
class]{Modular vector f\/ields and modular
class}\label{modularpoisson}

We consider a Poisson manifold $(M,\pi)$.
Recall that we denoted the Hamiltonian vector f\/ield with Hamiltonian
$f$ by $H^\pi_f$,
that the Poisson cohomology operator is $\d_{\pi} = [\pi , \cdot
\,]_{M}$ and that the Poisson homology operator is $\partial_\pi = [\d,
i_\pi]$.

Assuming that $M$ is
orientable, we choose a
volume form, $\lambda$, on $M$. The {\it{divergence}},
$\text{div}_{\lambda} Y$, of a
vector f\/ield $Y$ with respect to $\lambda$ is def\/ined by
$\mathcal{L}_Y \lambda =
(\text{div}_{\lambda} Y) \lambda$. Let us consider the linear map,
\begin{equation}\label{hamilt}
{X_\lambda : f  \mapsto \text{div}_{\lambda} (H^{\pi}_{f})} .
\end{equation}
Then
\begin{proposition}
\noindent$\bullet$ $X_\lambda$ is a derivation of
$C^{\infty}(M)$, i.e., a vector field on $M$,

$\bullet$ $X_\lambda$ is a $1$-cocycle in the Poisson
cohomology of $(M,\pi)$,

$\bullet$
the Poisson cohomology class of $X_\lambda$ is
independent of the choice of volume form, $\lambda$.
\end{proposition}

\begin{definition}
The vector f\/ield $X_\lambda$ is called a
\emph{modular vector field} of $(M,\pi)$.
The $\d_{\pi}$-cohomology class of $X_\lambda$ is called the
{\emph{modular class of the Poisson manifold}}
$(M,\pi)$.
\end{definition}

We shall denote the modular class of $(M,\pi)$ by $\theta(\pi)$.
The property $\d_\pi X_\lambda = 0$ means that $\mathcal L_{X_{\lambda}} \pi = 0$,
i.e., that the modular vector f\/ields are inf\/initesimal
automorphisms of the Poisson structure.

When $M$ is not orientable, densities are used instead of volume forms
in order to def\/ine the modular class.

\begin{example}\label{example1-sec2} If $(M,\pi)$ is symplectic, then $\theta(\pi) = 0$.
In fact, the Liouville volume form is invariant under all Hamiltonian
vector f\/ields.
\end{example}

If $M = \mathfrak{g}^{*}$, where $\g$ is a Lie
algebra, is equipped with the linear Poisson structure
(Kirillov--Kostant--Souriau Poisson structure),
the modular vector f\/ield associated with the standard
Lebesgue measure, $\lambda$, on the vector space, $\g^*$, is the constant vector
f\/ield on $\g^*$, i.e., linear form on $\g$, $X_\lambda = {\rm
  Tr}({\rm {ad}})$. Thus the modular vector f\/ield of the Poisson
manifold $\g^*$ is the
{\it {infinitesimal modular character}} of $\mathfrak{g}$.

\begin{example}\label{example2-sec2}
Consider the $2$-dimensional non-Abelian Lie algebra $\g$, with
basis $(e_1, e_2)$ and commutation relation, $[e_1,e_2]=e_2$. The trace
of the adjoint representation of $\g$ is the $1$-form  on
$\g$, $ e_1^*$.
In the dual vector space $\g^* \simeq {\mathbb R}^2$, with
dual basis $(u_1, u_2)$, where $u_1 = e_1^*$, $u_2 = e_2^*$,
and coordinates $(\xi_1 , \xi_2)$, the linear Poisson structure is def\/ined
by the bivector $\pi$ such that $\pi_{(\xi_1,\xi_2)} = \xi_2 \, u_1 \wedge
u_2$. The constant vector f\/ield $ u_1$ is the modular vector f\/ield of
$(\g^*,\pi)$ with
respect to the measure $\d \xi_1 \wedge \d \xi_2$. Since
$u_1$ is not globally Hamiltonian, the modular class of
$({\mathbb R}^2, \pi)$ is non-vanishing.
\end{example}

We have thus given an example of a Poisson structure
on ${\mathbb R}^2$ whose modular class is non-vanishing.
See~\cite{Radko1,Radko2} for examples of Poisson structures
with non-vanishing modular class on surfaces.

\subsection[Properties of modular vector fields]{Properties of modular vector f\/ields}

We state without proof the main properties of the modular vector
f\/ields of Poisson manifolds (see
\cite{Weinstein,Xu1,yks2000}). Each of the properties listed below
can be adopted as a def\/inition of the modular vector f\/ield,
$X_\lambda$.

\noindent
$\bullet$
For all $\alpha \in \Gamma (T^{*}M)$,
\begin{equation}\label{I}
{\langle \alpha , X_\lambda \rangle \lambda
= \mathcal{L}_{\pi^{\sharp}\alpha} \lambda
- (\iota_\pi \d  \alpha)\lambda}  ,
\end{equation}
a relation which reduces to \eqref{hamilt},
$X_\lambda (f) = {\rm{div}}_\lambda (H^\pi_f)  $, when $\alpha =
\d f$.

\noindent
$\bullet$
Let $n$ be the dimension of $M$. Given a
volume form, $\lambda$, the isomorphism
\[
*_\lambda :
\Gamma (\wedge^{\bullet} TM) \to \Gamma(\wedge^{n-\bullet}
T^{*}M)
\]
is def\/ined by
$*_{\lambda} u = \iota_{u} \lambda$,
for $u \in \Gamma (\w ^\bullet TM)$.

The modular vector f\/ield $X_\lambda$ is related to the $(n-1)$-form,
$\partial_\pi \lambda$, by
\begin{equation}\label{II}
{*_{\lambda} \, X_\lambda = - \partial_{\pi} \lambda} .
\end{equation}

\noindent
$\bullet$
Let us consider the operator on
$\Gamma (\w ^\bullet TM)$, of degree $-1$,
\[
\partial_{\lambda} = - (*_{\lambda})^{-1} \circ \d \circ
*_{\lambda} .
\]
On vector f\/ields, this operator coincides with
$- {\rm{div}}_\lambda$.
Applied to the Poisson bivector, it yields the modular vector f\/ield,
\begin{equation}\label{III}
{X_\lambda = \partial_{\lambda}\pi} .
\end{equation}
Thus the modular vector f\/ield $X_\lambda$ can be considered as the
``divergence'' of the Poisson bivector~$\pi$.

\noindent
$\bullet$
Let us also consider the operator
on $\Gamma (\w ^\bullet T^*M)$, of degree $-1$,
\[
  \partial_{\pi , \lambda} = - *_{\lambda} \circ \, \d _{\pi} \circ
(*_{\lambda})^{-1}  .
\]
Both operators
$\partial_{\pi}$ and $\partial_{\pi,\lambda}$ are {\it generators of
  square zero} of the Gerstenhaber algebra \linebreak
$(\Gamma (
\wedge^{\bullet} T^{*}M) , [~,~]_{\pi})$. In view of equation~\eqref{generator},
the dif\/ference
$\partial_{\pi,\lambda}- \partial_{\pi}$ is the
interior pro\-duct by a vector f\/ield. In fact, this vector f\/ield
is a $1$-cocycle in the
$\d_\pi$-cohomology \cite{Koszul,Xu1,yks2000}, and it coincides with the
modular vector f\/ield,
\begin{equation}\label{IV}
{\partial_{\pi,\lambda}- \partial_{\pi} = \iota_{X_\lambda}} .
\end{equation}

\section{The modular class of a
Lie algebroid with a Poisson structure}\label{LAwithpoisson}
It is straightforward to generalize the def\/inition of the modular
class from the case of the Poisson manifolds to that of the
Lie algebroids with a Poisson structure. Let
$(A,a_A,[\ ,\ ]_{A})$ be a Lie algebroid
with a Poisson bivector
$\pi \in \Gamma (\wedge^{2} A)$.
Since the exact $1$-forms,
$\d _A f$, $f \in C^\infty (M)$, do not, in general, span the space of
sections of $A^*$, we cannot rely on the original def\/inition
\eqref{hamilt}, but, assuming that $A$ is orientable,
we choose a nowhere-vanishing section $\lambda$ of
$\w^{\rm{top}}(A^*)$, and we consider the expressions in formulas \eqref{I}--\eqref{IV},
where $\d$ is replaced by $\d _A$, $\mathcal L$ by
$\mathcal L^A$, where $\mathcal L^A_X = [\d_A, \iota_X]$ for $X \in \Gamma
A$, $\partial_\lambda$
by $- (*_\lambda)^{-1} \circ \d_A \circ *_\lambda $,
etc. Each of the equations~\eqref{I}--\eqref{IV}
def\/ines uniquely the same section $X_\lambda$ of
$A$ which is a $1$-cocycle in the $\d_\pi$-cohomology. Furthermore,
the class of $X_\lambda$ is
independent of the choice of section of $\w^{\rm{top}}(A^*)$.
The case of a~Poisson manifold~$M$ is recovered when $A = TM$.

\begin{definition}\label{modclasspoisson}
The section $X_\lambda$ of $A$ is called a
\emph{modular section of the Lie algebroid with a Poisson
  structure $(A,\pi)$}.
The $\d_{\pi}$-cohomology class of $X_\lambda$
is called the {\emph{modular class}} of $(A,\pi)$.
\end{definition}

We shall denote the modular class of $(A,\pi)$ by
$\theta(A,\pi)$. When $A=TM$, the modular class of $(A,\pi)$ is the
modular class $\theta(\pi)$ of the Poisson manifold $(M,\pi)$.

The non-orientable case can be dealt with by using densities instead
of volume forms.

In the next section, we shall def\/ine the modular class,
${\rm {Mod}} \, E$,
of a Lie algebroid $E$, and in Section \ref{revisited} we shall show
how $\theta(A,\pi)$ is related to ${\rm {Mod}} (A^*)$ when
$A^*$ is the Lie algebroid with anchor $a_A \circ \pi^\sharp$ and bracket
$[~,~]_\pi$.

\section{The modular class of a Lie algebroid}\label{section4}
The modular class of a Lie algebroid was introduced by Weinstein in
\cite{Weinstein}.
This section summarizes some of the results of \cite{ELW}, in which
the theory was developed.

\subsection{Lie algebroid representations}

A {\it {representation}} \cite{M1,M2} of a Lie algebroid
$(E, a_E, [~,~]_E)$ with base $M$
in a vector bundle $V$ with base $M$ is a map $D : \Gamma E \times \Gamma V
\to \Gamma V$, $(u,s) \mapsto D_{u} s$, such that
\begin{gather*}
D_{fu} s  = f D_{u} s , \\
D_{u} (fs)  = fD_{u} s+ (a_E(u) \cdot f)s  , \\
D_{[u,v]_{E}}  = [D_{u},D_{v}]  ,
\end{gather*}
for all $u$, $v \in \Gamma E$, $s \in \Gamma V$ and $f \in C^\infty(M)$.
A representation of $E$ in $V$ is also called a {\it flat
  $E$-connection}
\cite{Xu1}  on
$V$, and $V$ is called an {\it $E$-module} \cite{Vaintrob, Fernandes1} if there is a representation of~$E$ on $V$. Clearly, if $E = TM$, a f\/lat $E$-connection is
a f\/lat connection in the usual sense,
and if
$M = \{{\rm {pt}}\}$, this
notion reduces to the Lie
algebra representations in vector spaces.

The {\it canonical representation} \cite{ELW} of a Lie algebroid $E$ is
the representation $D^E$ of $E$ in the line bundle
$L^E = \wedge^{\rm top}E \otimes
\wedge^{\rm top} (T^{*} M)$ def\/ined by
\begin{equation}\label{canon}
D^{E}_{u} (\lambda \otimes \nu) = \mathcal{L}^{E}_{u} \lambda \otimes
\nu + \lambda \otimes \mathcal{L}_{a_E(u)} \nu \ ,
\end{equation}
where $u \in \Gamma E$, $\lambda \in \Gamma
(\wedge^{\rm top}E)$, $\nu \in \Gamma
(\wedge^{\rm top} T^{*} M)$,
and where ${\mathcal L}^E_{u} \lambda = [u ,\lambda]_{E}$.

When $E$ is the integrable sub-bundle of $TM$ def\/ining a foliation
$\mathcal F$ of
$M$, the line bundle $L^E$ is isomorphic to the top exterior power of
the conormal bundle to $\mathcal F$, and the representation $D^E$ is
the top exterior power of the Bott connection of $\mathcal F$.
When the normal bundle is
orientable, the sections of $L^E$ are the tranverse measures to the
foliation.

\subsection{Characteristic
class of a Lie algebroid representation in
a line bundle}

Let $E$ be a Lie algebroid
with a representation $D$ in
a line bundle $L$, and let $s$ be a
nowhere-vanishing section of $L$.
Def\/ine $\theta_{s} \in \Gamma (E^{*})$
by the condition
\[
{ \langle \theta_{s} , u \rangle s =
  D_{u} s}  ,
\]
for all $u \in \Gamma E$.
Then $\theta_{s}$ is a $\d_{E}$-{\it cocycle}. Furthermore, the class
of $\theta_s$ is
independent of the choice of section $s$ of $L$, whence the following
def\/inition \cite{ELW}.

\begin{definition}
The section
$\theta_s$ of $E^*$ is
called a {\it characteristic cocycle} associated with the
  representation $D$ and the section $s$. Its class is called the {\it
  characteristic class} of the representation~$D$.
\end{definition}

If $L$ is not trivial, the class of $D$ is def\/ined as one-half that of the
class of the associated representation of $E$ in the square of $L$.

\subsection{The modular class of a Lie algebroid}\label{4.3}
Equation \eqref{canon} def\/ines the canonical representation of the Lie
algebroid $E$ in the line bundle, $L^E =
\wedge^{\rm top}E \otimes
\wedge^{\rm top} T^{*} M$, which appears in the following def\/inition \cite{ELW}.

\begin{definition}
The characteristic class of the canonical representation $D^E$ of $E$
is called the {\it modular class} of the Lie algebroid $E$.
A cocycle belonging to the modular class of $E$ is callled a~{\it{modular cocycle}} of $E$.
\end{definition}

We shall
denote the modular class of $E$ by ${\rm Mod} E$.
Thus, by def\/inition, in the orientable case,
$\text{Mod } E$ is the $\d_{E}$-class of the 1-cocycle $\theta
\in \Gamma (E^{*})$, depending on the nowhere-vanishing section
$\lambda \otimes \nu$ of $L^E$, such that
\[
\langle \theta , u \rangle \lambda \otimes
\nu = \mathcal{L}^{E}_{u} \lambda \otimes \nu + \lambda \otimes
\mathcal{L}_{a_E(u)}  \nu ,
\]
for all $u \in \Gamma E$.

\begin{example}\label{example1-sec4} If $E = TM$, the
modular class vanishes,
${\rm Mod}\, (TM) = 0$, since one can choose $\lambda \in \Gamma
(\w^{\rm{top}}T^*M)$ and $\nu \in \Gamma
(\w^{\rm{top}}TM)$ such that $\langle \lambda, \nu \rangle = 1$.
\end{example}

If $E$ is a Lie algebra $\g$, then
$\text{Mod }\mathfrak{g}$ is the {\it{infinitesimal modular
character}},
$\rm{Tr} ({\rm {ad}})$,  of $\g$.
Comparing this fact with the result
in Section \ref{modularpoisson},
we see that
the modular class of $\g$,
considered as a Lie algebroid over a point, is equal to the modular
class of the Poisson manifold~$\g^*$.

\begin{example}\label{example2-sec4}
If $\g$ is the non-Abelian 2-dimensional Lie algebra
with commutation relation $[e_1,e_2] = e_2$, then $\text{Mod} \, \mathfrak{g}
= e_1^*$. We saw in Example~\ref{example2-sec2} of Section~\ref{modularpoisson} that $e_1^*$ is
a modular vector f\/ield of the Poisson manifold~$\g^*$.
\end{example}

The preceding property of Lie algebras
can be extended to Lie algebroids in the following way.
Let us identify each section of a vector bundle $F \to M$ with its
vertical lift,
which is a~section of the vertical tangent bundle $VF \to F$. This
vertical lift, in turn, can be identif\/ied with a vector f\/ield on $F$,
tangent to the f\/ibers and invariant by translations along each f\/iber.
In this way, the modular class of a Lie algebroid $E$ can be
identif\/ied~\cite{Weinstein, NW} with
the modular class of the Poisson manifold $E^*$, when $E^*$ is equipped
with the {\it{linear Poisson structure}} (see, e.g.,~\cite{M2})
def\/ined by the anchor and bracket of $E$.

\begin{example}\label{example3-sec4} When $E$ is the Lie algebroid associated to a
foliation $\mathcal F$, and the normal bundle of $\mathcal F$ is
orientable, the sections of $L^E$ are the tranverse measures to the
foliation. Therefore the modular class of $E$
is the obstruction to the existence of a transverse measure invariant
under all vector f\/ields tangent to the foliation.
\end{example}

\begin{remark}
\label{remark-sec4}
A section of a vector bundle $E^*$ can be considered as a function on
$E$ viewed as a~supermanifold, and a nowhere-vanishing section of $L^E$ determines a
Berezinian
volume on the supermanifold $E$.
A modular cocycle of a Lie algebroid $E$ is
the divergence with respect to such a Berezinian volume of $\d_E$
considered as a vector f\/ield on the supermanifold $E$ (see~\cite{ELW}).
\end{remark}

\subsection{The case of Poisson manifolds}\label{pmanifolds}

If $(M,\pi)$ is a Poisson manifold, the comparison of the modular class,
$\theta(\pi)$, of $(M,\pi)$ with the modular class, $\text{Mod }
(T^{*}M)$, of the Lie algebroid
$(T^{*} M , \pi^{\sharp} , [~,~]_{\pi})$ yields the following
result which was proved in \cite{ELW}.

\begin{theorem}\label{theorem1}
For a Poisson manifold, $(M,\pi)$, the modular classes
$\theta(\pi)$ and ${\rm{Mod}} \, (T^{*}M)$ are equal, up to a factor $2$,
\begin{equation}\label{poissononehalf}
{\theta(\pi) = \frac{1}{2 } \, {\rm{Mod}} \, (T^{*}M) } .
\end{equation}
\end{theorem}

If $A$ is a Lie algebroid with a Poisson structure def\/ined by $\pi \in
\Gamma (\w^2 A)$, the question arises whether the
cohomology class, $\theta(A ,\pi)$, that was def\/ined in
Section \ref{LAwithpoisson}, and $\rm{Mod}(A^{*})$, the modular class
of the Lie algebroid
$(A^{*}, a_A \circ \pi^{\sharp} , [~,~]_{\pi})$
satisfy a relation as simple as \eqref{poissononehalf}.
The answer is no in general, and the correct relation is obtained
by the introduction of a new notion, the modular class of a morphism.

\section{The modular class of a Lie algebroid
  morphism}\label{section5}

This section deals f\/irst with the
modular classes of base-preserving Lie algebroid morphisms which
were introduced in
\cite{GMM} and \cite{yksW}.
In Section
\ref{general}, we shall show how the def\/inition and properties of
these classes can be extended to general morphisms.

\subsection{Lie algebroid morphisms}

By def\/inition, if $E$ and $F$
are Lie algebroids over the same base, a vector bundle map
$\Phi : E \to F$ is a {\it Lie algebroid
morphism} (over the
identity of $M$)
if $a_E = a_F \circ \Phi$ and $\Phi$ induces a Lie algebra
homomorphism from $\Gamma E$ to $\Gamma F$.
It is well known  (see, e.g., \cite{yksM,Huebschmann1})
that $\Phi$ is a morphism if and only if
$\wedge^{\bullet}(\Phi^{*})$ def\/ines
a {\it chain map} from
$(\Gamma (\wedge^{\bullet} F^{*}) , \d_{F})$ to
$(\Gamma(\wedge^{\bullet}E^{*}) , \d_{E})$.

Given a morphism $\Phi : E \to F$, let us consider
the well-def\/ined class in the Lie
algebroid cohomology of $E$,
\[
{\text{Mod }({\Phi}) = \text{Mod } E - \Phi^{*} (\text{Mod } F )} .
\]
\begin{definition}\label{def}
The $\d_E$-cohomology class $\text{Mod }({\Phi})$ is called the
{\emph{modular class of the Lie algebroid morphism}} $\Phi$.
\end{definition}

The modular class of $\Phi$ can be
considered as a {\it relative modular class} of the pair
$(E,F)$, whence the terminology adopted in \cite{yksW}.
A representative of the class $\text{Mod }({\Phi})$ is a section of~$E^*$. It is clear that the modular class of an isomorphism of Lie
algebroids vanishes.

It was proved in \cite{yksW} that the modular class of a morphism $\Phi : E \to F$ is the
characteristic class of a representation of $E$.
In fact, set $L^{E,F} = \wedge^{\rm {top}} E \otimes \wedge^{\rm
  {top}} F^*$ and def\/ine the map $D^{\Phi}$ by
\begin{equation}\label{repmorphism}
D^{\Phi}_{u} (\lambda\otimes \nu) = \mathcal{L}^{E}_{u}\lambda
\otimes \nu
+ \lambda \otimes \mathcal{L}^{F}_{\Phi u} \nu ,
\end{equation}
for $u\in \Gamma E$ and
$\lambda \otimes \nu \in \Gamma(L^{E,F})$.
\begin{proposition}\label{charmorphism}
When $\Phi : E \to F$ is a Lie algebroid morphism,
the map $D^\Phi$ is a representation of $E$ on
the line bundle $L^{E,F} $, and the modular class of
$\Phi$ is the characteristic class of this representation.
\end{proposition}

\subsection{The modular class of a Lie algebroid with a Poisson
  structure revisited}\label{revisited}

Let us again consider
a Lie algebroid $A$ with a Poisson structure def\/ined by $\pi \in
\Gamma (\w^2 A)$, as in Section \ref{LAwithpoisson}. We can now
establish a relation between
the cohomology class $\theta(A ,\pi)$ and $\rm{Mod}(A^{*})$, the modular class
of the Lie algebroid
$(A^{*}, a_A \circ \pi^{\sharp} , [~,~]_{\pi})$.
(The modular class $\theta(A,\pi)$ is the class of a section of
$A$, not to be confused with the modular class of $A$ as a Lie
algebroid, ${\rm{Mod}}\, A$, which is the class of a section of $A^*$.)

\begin{theorem}\label{untwisted}
The modular class of a Lie algebroid with Poisson structure,
$(A,\pi)$, and the modular class of the morphism $\pi^\sharp:
A^* \to A$ are equal up to a factor $2$,
\begin{equation} \label{onehalf}
{\theta(A,\pi) =\frac{1}{2} \, \rm{Mod }({\pi^\sharp})} .
\end{equation}
\end{theorem}

The equality
\eqref{onehalf}
yields the desired relation between $\theta(A ,\pi)$
and $\rm{Mod}(A^{*})$,
\begin{equation}\label{V}
\theta(A,\pi) = \frac{1}{2} \big(\text{Mod } (A^{*}) -
  (\pi^{\sharp} )^{*} (\text{Mod } A) \big) .
\end{equation}
In particular, if $\pi$ is non-degenerate, $\theta(A, \pi)$ vanishes.

In the case of a Poisson manifold, $A=TM$. Since $\text{Mod } (TM) = 0$,
\eqref{V} reduces in this case to
$\theta(\pi) = \frac{1}{2} \, \text{Mod } (T^*M)$, i.e.,
we recover~\eqref{poissononehalf}, the result of Theorem~\ref{theorem1}.

If $A = \g$ is a Lie algebra
equipped with a solution $\pi = r$ of the classical Yang--Baxter
equation~\eqref{cybe}, then $(\g,r)$ can be considered as a Lie
algebroid over a point with a Poisson structure, and as such it has
a modular class, $\theta(\g,r)$.

\begin{example}\label{example-sec5} Let $\g$ be the $2$-dimensional non-Abelian Lie
algebra as in Example~\ref{example2-sec2} of Section~\ref{modularpoisson}. Then $r = e_1 \wedge
e_2$ is a solution of the classical Yang--Baxter equation on $\g$. Because $r$
is non-degenerate, $\theta(\g,r) =0$. On the other hand, the modular class of the
Lie algebra $\g^*$ equipped with the Lie bracket $[~,~]_r$ def\/ined by
\eqref{rmatrix} can be computed as the trace of the adjoint representation
of this Lie algebra, yielding $\text{Mod } (\g^{*}) = - e_2$.
Since $\text{Mod } \g = e_1^*$, we again obtain
${\theta(\g,r) =
\frac{1}{2} \left (\text{Mod } (\g^{*}) +
  r^{\sharp}  (\text{Mod } \g) \right)} = \frac{1}{2} (- e_2 +
r^\sharp e_1^*) = 0$.
\end{example}

\subsection{Unimodularity}

A Lie algebroid or a Lie algebroid morphism
is called {\it{unimodular}} if its modular class
vanishes. Examples of unimodular Lie algebroids
are tangent bundles and unimodular Lie algebras, and we have
mentioned above in Example~\ref{example3-sec4} of Section \ref{4.3}
the meaning of unimodularity in the theory of foliations.
The unimodularity of morphisms also is related to the existence of
invariant transverse measures. For example,
if $H$ is a connected, closed subgroup of a connected Lie group
$G$, with Lie algebras $\mathfrak h \subset \g$,
the canonical injection
$i : \mathfrak{h} \hookrightarrow \mathfrak{g}$ is unimodular if
and only if there exists a~$G$-invariant measure on the homogeneous
space $G/H$ \cite{yksW}.

\begin{remark}\label{remark1-sec5} When viewed as a supermanifold, an orientable
unimodular Lie algebroid is a~QS-manifold in the sense of Schwarz
\cite{Schwarz} (see \cite{yksMo}).
\end{remark}

A Poisson structure, $\pi$, on a Lie algebroid, $A$, is also called
{\it{unimodular}} if the class $\theta(A,\pi)$ vanishes.
A Poisson manifold $(M,\pi)$ is {\it{unimodular}}
if $\theta(\pi) = 0$, i.e.,
if there exists a
vanishing modular vector f\/ield $X_\lambda$,
i.e., if there exists a
density that is invariant under all Hamiltonian vector f\/ields.

On an orientable
Lie algebroid with a Poisson structure, the Poisson cohomology
$H^\bullet(A,\d_\pi)$ and the Poisson homology
$H_\bullet(A,\partial_\pi)$ can be compared.
When $\lambda$ is a volume element on $A$,
$*_\lambda$ is a chain map from the complex
$(\Gamma(\w^\bullet A^*),\partial_\pi)$ to the complex
$(\Gamma(\w^{n-\bullet} A),\d_\pi + \epsilon_{X_\lambda})$, where
$n$ is the rank of $A$ and $\epsilon_X$ is the exterior
product of multivectors by the section $X$ of $A$.
In fact, in view of equation~\eqref{IV}, $\partial_{\pi} =
\partial_{\pi, \lambda} -
\iota_{X_\lambda} = - *_{\lambda} \circ \, (\d_{\pi} +
\epsilon_{X_\lambda}) \circ
(*_{\lambda})^{-1}$.
Therefore, the homology of the complex
$(\Gamma(\wedge^{\bullet}A^{*}),\partial_{\pi})$ is isomorphic to the
cohomology of the complex $(\Gamma(\wedge^{\bullet} A), \d_{\pi} +
\epsilon_{X_\lambda})$
\cite{Huebschmann3,Huebschmann2,Xu1,yks2000}.
In particular,
\begin{proposition}
When the Poisson structure, $\pi$, of an orientable Lie algebroid,
$A$, is {\it{unimodular}},
its Poisson homology and Poisson cohomology are isomorphic,
\[
H_{\bullet}(A,\partial_{\pi}) \simeq H^{{\rm {\rm top}}- \bullet}
(A, \d_{\pi}) .
\]
\end{proposition}

\begin{remark}\label{remark2-sec5} More general pairings between Lie algebroid
cohomology and
homology are to be found in \cite{Huebschmann2,ELW,Xu1}.
In \cite{GMM}, Grabowski, Marmo and Michor use odd volume forms
to def\/ine generating operators of square zero
of the Gerstenhaber algebra $\Gamma (\w^\bullet A)$ for an arbitrary Lie
algebroid $A$, and they show that the associated homology
is independent of the choice of odd volume form and that, when
the Lie algebroid is orientable, there is a Poincar\'e duality
between this homology and the Lie algebroid
cohomology.
\end{remark}

\subsection{General morphisms of Lie algebroids}\label{general}
The extension of the def\/inition and properties of the modular class
to the case of Lie algebroid morphisms that are not necessarily
base-preserving is the subject of the article \cite{yksLW}.
Let us
brief\/ly review this general case.
When $\Phi$ is a vector bundle map from $E \to M$ to $F \to
N$ over a map $\phi : M \to N$, there is a map $\widetilde\Phi^*$ from
the
sections of $F^* \to N$ to the sections of $E^* \to M$, induced from the
base-preserving vector bundle morphism $\Phi^* : \phi^! F^* \to E^*$,
where $\phi^! F^*$ is the pull-back of $F^*$ under $\phi$.
Then $\widetilde\Phi^*$ is extended as an exterior
algebra homomorphism $\w^\bullet \widetilde\Phi^* : \Gamma(\w^\bullet
F^*) \to \Gamma(\w^\bullet E^*)$. Let $E \to M$ and $F \to N$ be Lie
algebroids. By def\/inition, $\Phi$ is a~{\it Lie algebroid morphism} if
$\wedge^{\bullet} \widetilde{\Phi}^{*}$ is
a {\it {chain map}} from
$(\Gamma (\wedge^{\bullet} F^{*}) , \d_{F})$ to
$(\Gamma(\wedge^{\bullet}E^{*}) , \d_{E})$.

\begin{remark}\label{remark3-sec5}
 It was proved by Chen and Liu \cite{CL} that this def\/inition is
equivalent to the original def\/inition of morphism
in \cite{HM,M1,M2}. In fact, it is also
equivalent to the def\/inition in \cite{Vaintrob}: $\Phi$~maps
the homological vector f\/ield $\d_E$ on the supermanifold $E$ to $\d_F$.
\end{remark}

It follows from the chain map property that the $\d_E$-cohomology class,
\[
{\rm Mod} (\Phi) = {\rm Mod} \, E - \widetilde \Phi^* ({\rm Mod} \, F) ,
\]
is well-def\/ined. It is called the
{\it modular class of the morphism} $\Phi$, a def\/inition that
generalizes Def\/inition \ref{def}.

That this class is the characteristic class of a
representation of $E$ is proved using the non-trivial fact that the
pull-back of a Lie algebroid
representation by a morphism exists.
Another theorem in \cite{yksLW} states that the pull-back of a
Lie algebroid by a transverse map, in particular by a submersion,
gives rise to a
morphism with vanishing modular class.

\section{The regular case and triangular Lie
  bialgebras}\label{regulartriang}
A Lie algebroid $A$ over $M$
is called {\it regular} if
the image of its anchor map
is of constant rank, and thus a sub-bundle of $TM$.
The Lie algebroids associated to
foliations are an important example of regular Lie algebroids.
The modular classes of regular Lie algebroids were considered in
\cite{Crainic, Fernandes0, Fernandes1}.

A Poisson manifold $(M,\pi)$ is called {\it regular}
if its cotangent Lie algebroid $T^*M$ is regular, i.e.,
if the image of the map $\pi^\sharp : T^*M \to TM$ is of constant
rank. By extension, a Poisson structure,~$\pi$, on a Lie algebroid,
$A$, is called {\it
regular} if the image of $\pi^\sharp$ is of constant rank, and
thus a Lie subalgebroid of $A$.
We now summarize results on the modular class of Lie algebroids with a
regular Poisson structure obtained in \cite{yksY}.

\subsection{The regular case}\label{sectionregular}
Let
$(A ,\pi)$ be a Lie algebroid with a
regular Poisson structure.
Let $B = \pi^{\sharp} (A^{*})$ and let us denote the Lie algebroid morphism
$\pi^\sharp$ considered as a submersion onto $B$
by $\pi^{\sharp}_{B}$.
There is an exact sequence of Lie algebroids
over the same base,
\[
0 \to {\rm{Ker}} \, \pi^{\sharp}
\to A^{*} \overset{\pi^{\sharp}_{B}}{\longrightarrow} B \to 0 .
\]
The {\it canonical representation} of $B$ on ${\rm{Ker}} \, \pi^{\sharp}$ is obtained by
factoring the
adjoint action of $\Gamma(A^{*})$, with respect to the Lie bracket
$[~,~]_\pi$ def\/ined by \eqref{FMMbracket},
on $\Gamma({\rm{Ker}} \, \pi^{\sharp})$ through the submersion
$\pi^{\sharp}_{B} : A^{*} \to B$. (The existence of this
representation follows from \cite{M2}, Proposition 3.3.20.)
Explicitly,
\[
 X \cdot \gamma = \mathcal{L}^{A}_{X} \gamma ,
\]
for $X \in \Gamma B$, $\gamma \in \Gamma ({\rm{Ker}} \, \pi^{\sharp})$.
\begin{theorem}\label{theoremregular}
The modular class, $\theta(A,\pi)$, of a Lie algebroid with a regular
Poisson structure, $(A,\pi)$, satisf\/ies
\begin{equation}\label{regular}
{\theta(A , \pi) = (\pi^{\sharp}_{B})^{*} (\theta_{B})}  ,
\end{equation}
where $\theta_{B}$ is the
characteristic class of the representation of $B$ in
$\wedge^{\rm top} ({\rm{Ker}} \, \pi^{\sharp})$ induced from the canonical representation.
\end{theorem}
This theorem
permits an ef\/f\/icient computation of a representative of the
modular class, avoiding the computation of terms which
mutually cancel in equation \eqref{V}.

\begin{remark}\label{remark1-sec6}
Formula \eqref{regular} can also be obtained from the more general
computation of the modular class of a Lie algebroid
morphism with constant rank and unimodular kernel in~\cite{yksLW}.
\end{remark}

\subsection{Triangular Lie bialgebras}
This section summarizes some of the results of \cite{yksY} concerning
f\/inite-dimensional real or complex Lie
algebras with a triangular $r$-matrix, the def\/inition of which was
recalled in Section \ref{triangular}.

Let $\g$ be a Lie algebra considered as a Lie algebroid over a
point, and let $r \in \w^2 \g$ be a~triangular $r$-matrix on $\g$,
def\/ining a triangular Lie bialgebra structure on $(\g,\g^*)$.
The image of $r^\sharp:
\g^* \to \g$
is called
the {\it{carrier}}
of $r$. We set $\p = {\rm {Im}} \, r^\sharp$.
In this case, the canonical representation of $\p$ on ${\rm
  {Ker}} \, r^\sharp$ coincides with the restriction of the coadjoint
representation of $\g$, and it is dual to the representation of
$\p$ on $\g/\p$
induced from the adjoint action of $\g$. Therefore, as a conseqeunce of
Theorem  \ref{theoremregular},
we obtain the following proposition, where
$\chi_{\mathfrak{p},V}$ denotes the {\it infinitesimal character} of a
representation $V$ of Lie algebra $\mathfrak{p}$, and
${r^{\sharp}_{(\p)}}$ denotes $r^\sharp$ considered as a
skew-symmetric isomorphism from
$\p^{*} = \g^{*}/{\rm {Ker}} \, r^{\sharp}$ to $\p$.
\begin{proposition}\label{bialgebra}
Let $r$ be a triangular $r$-matrix on a Lie algebra $\g$. The modular
class of $(\g,r)$ is the element $\theta(\g,r)$
of $\g$ such that
\[
{\theta (\mathfrak{g}, r ) = -  r^{\sharp}_{({\mathfrak p})}
  (\chi_{{\mathfrak p} , {\rm{Ker}} \, ( r^{\sharp})}) =
  r^{\sharp}_{({\mathfrak p})} (\chi_{{\mathfrak p},
    \mathfrak{g}/{\mathfrak p}})} .
\]
\end{proposition}

\subsection{Frobenius Lie algebras}\label{frobalg}
The following proposition was proved by Stolin \cite{Stolin} (see also
\cite{GG,HY}).
\begin{proposition}\label{stolin}
Let $\omega \in
\wedge^2 (\mathfrak p^*)$ be a
non-degenerate $2$-form on a Lie algebra $\p$, and let $r \in \w^2 \p$
be the inverse of $\omega$.
Then
${\rm d}_\p \omega = 0$ if and only if $[r,r]_\g = 0$.
\end{proposition}

A Lie algebra
$\mathfrak p$ is called a \emph{quasi-Frobenius} (resp.,
\emph{Frobenius}) \emph{Lie algebra} if there exists a~non-degenerate
$2$-cocyle (resp., $2$-coboundary) on $\mathfrak p$.
A Lie algebra
$\mathfrak p$ is called \emph{Frobenius with respect to} $\xi \in
\mathfrak p^*$ if $\omega = - {\rm d}_\mathfrak p \xi$ is a
non-degenerate $2$-form.

From this def\/inition and from Proposition~\ref{stolin}, we obtain:
\begin{proposition}
Let $\g$ be a Lie algebra, and let $\p$ be a Lie subalgebra of $\g$.
Assume that $\mathfrak p$ is Frobenius with respect to $\xi \in
\mathfrak p^*$. Then
$r = - ({\rm d}_\mathfrak p \xi)^{-1} \in \w^2 \p$
is a triangular $r$-matrix on $\g$, i.e., a
solution of the classical Yang--Baxter equation.
\end{proposition}
 When the carrier of a triangular $r$-matrix is a Frobenius Lie
 algebra, Proposition \ref{bialgebra} yields the following result.
\begin{theorem}\label{frobenius}
Let $r$ be a triangular $r$-matrix on a Lie algebra $\g$. Assume that
$\p =  {\mathrm{Im}} \,  r^\sharp$ is Frobenius with respect to $\xi
\in \p^*$. Then $\theta(\g, r)$ is the unique element
$X \in \g$ such that
\begin{equation}\label{character}
{\rm{ad}}^*_X \, \xi = \chi_{{\mathfrak p}, \mathfrak{g}/{\mathfrak
    p}} .
\end{equation}
\end{theorem}

\begin{example}\label{example1-sec6}
We illustrate this theorem with the example where $\p \subset
\sl_3(\mathbb R)$ is
the Lie algebra of traceless matrices of the form $\begin{pmatrix}\bullet & \bullet & \bullet \\
0 & \bullet & \bullet \\
0 & \bullet & \bullet
\end{pmatrix}$. It
is Frobenius with respect to
$\xi = e_{12}^* + e_{23}^*$.
The coboundary of $ - \xi$ is
$- \d_\p \xi = 2( e_{11}^* \w e_{12}^* +
e_{12}^* \w e_{22}^* + e_{13}^* \w e_{32}^* + e_{22}^* \w e_{23}^* +
e_{23}^* \w e_{33}^*)$,
and the corresponding triangular $r$-matrix is
\[
r_{GG} = (\frac{2}{3} e_{11}
-
\frac{1}{3} (e_{22} + e_{33})) \w e_{12}
+
\left(\frac{1}{3} (e_{11} + e_{22}) -
\frac{2}{3} e_{33} \right) \w e_{23} + e_{13} \w e_{32}  .
\]
The image of $r_{GG}^\sharp$ is $\p$ and its kernel
is generated by $e_{21}^*$ and $e_{31}^*$.
Computing the character
$\chi_{\p,\g/\p} = - 2 e^*_{11} + e^*_{22} + e^*_{33}$, and solving
equation~\eqref{character}, we obtain
$\theta(\sl_3(\mathbb R),r_{GG}) = -  2 e_{12} + e_{23}$.
The $r$-matrix $r_{GG}$ is the Gerstenhaber--Giaquinto generalized Jordanian
$r$-matrix \cite{GG} on $\sl_3(\mathbb R)$.
See \cite{yksY} for the computation of the modular class in the
general case of $\sl_n(\mathbb R)$.
\end{example}

\section{The modular class of a twisted Poisson structure}\label{section7}
Twisted Poisson structures on manifolds were studied by \v Severa
and Weinstein \cite{SW} and such structures on Lie algebroids were
subsequently
def\/ined by Roytenberg in \cite{Roytenberg}. It was proved there that a twisted
Poisson structure on $A$ corresponds to a
{\it quasi-Lie-bialgebroid}
structure on $(A,A^*)$,
generalizing the triangular Lie bialgebroids
def\/ined by Poisson structures.
The closed $3$-form which we denote by
$\psi$ below, to conform with the notation in \cite{yks1992} and
subsequent publications,
plays an important role in generalized complex geometry
\cite{Gualtieri}
and, consequently, in recent studies on sigma-models, where it is
usually denoted by $H$, and sometimes called an $H$-{\it flux}.

\subsection{Twisted Poisson structures}\label{7.1}
A {\it {twisted Poisson structure}} on a Lie algebroid $A$ is a pair
$(\pi,\psi)$, where $\pi \in \Gamma (\w^2 A)$ and $\psi$ is a~$d_A$-closed
form on $A$ satisfying
\begin{equation}\label{twisted}
\frac{1}{2} [\pi,\pi]_A = (\w^3\pi^\sharp) \psi .
\end{equation}
When a twisted Poisson structure is def\/ined on the Lie algebroid
$TM$, the manifold $M$ is called a {\it{twisted Poisson manifold}}.

\begin{example}\label{example1-sec7}
There is no genuine twisted Poisson structure on a Poisson manifold of
dimension $ \leq~3$, since the image of $\pi^\sharp$
is then of dimension at most $2$.
We give an example  \cite{Hirota} of a twisted Poisson
structure on
an open subset of $\mathbb R^4$,
$M = \{(x_1, x_2, x_3, x_4) \in \mathbb R^4 \, | \, x_1 \neq 0 \quad
{\mathrm {and}} \quad  x_3 \neq 0 \}$.
Denote $\frac{\partial}{\partial x_i}$
by $\partial_i$. Let $\pi = x_3 \partial_1 \wedge \partial_2 + x_1
\partial_3 \wedge  \partial_4$ and $\psi = ((x_1)^{-2}
{\rm d}x_4 - (x_3)^{-2} {\rm d}x_2) \wedge {\rm d}x_1
\wedge {\rm d}x_3$.
Then $(\pi, \psi)$ is a twisted Poisson structure on the manifold
$M$.
\end{example}

When $(A,\pi ,\psi)$ is a Lie algebroid with a
twisted Poisson structure, the
bracket def\/ined by
\[
[\alpha, \beta]_{\pi,\psi}
= \mathcal{L}^A_{\pi^{\sharp}\alpha}
\beta - \mathcal{L}^A_{\pi^{\sharp}\beta} \alpha - \d_A (\pi (\alpha
,\beta)) + \psi(\pi^\sharp \alpha , \pi^\sharp \beta , \cdot)  ,
\]
for $\alpha$ and $\beta \in \Gamma (A^*)$, satisf\/ies the Jacobi
identity. With the anchor $a_A \circ \pi^\sharp$ and the bracket
$[~,~]_{\pi,\psi}$,
$A^{*}$ is a Lie algebroid, and the map $\pi^{\sharp}$ is a Lie
algebroid morphism
from $A^{*}$ to $A$.

There is an associated
Gerstenhaber bracket on $\Gamma (\w^\bullet A^*)$, also denoted by
$[~,~]_{\pi,\psi}$,
and an associated dif\/ferential
 on $\Gamma(\w^\bullet A)$,
$\d_{\pi, \psi} = [\pi, \, \cdot \, ]_{A} + \underline{\d}$, where
$\underline{\d}$ is the interior product of sections of $\w^\bullet A$
by the bivector-valued $1$-form, $(\alpha, \beta) \mapsto
\psi(\pi^\sharp \alpha , \pi^\sharp \beta , \cdot)$.

\begin{remark}\label{remark1-sec7} The pair $(A,A^*)$ is not in general a Lie
bialgebroid. However, when $A$ is equipped with the bracket def\/ined by
$[X,Y]_{A,\psi} = [X,Y]_A + \psi(\pi^\sharp (\, \cdot \,) , X,Y)$,
for $X$ and $Y \in \Gamma A$, the pair $(A,A^*)$ becomes a
quasi-Lie-bialgebroid.
In terms of the big bracket \cite{Roytenberg,yks2007},
\begin{gather*}
[X,Y]_{A,\psi} = \{\{X,\mu_\psi\},Y\}  , \qquad {\rm{where}} \quad
\mu_\psi = \mu - \{\psi, \pi \}  , \\
[\alpha, \beta]_{\pi,\psi} = \{\{\alpha, \gamma_\psi\}, \beta\}  ,
\qquad {\rm{where}} \quad
\gamma_\psi = - \{ \mu ,\pi \} + \frac{1}{2} \{\{ \psi , \pi \}
,\pi \}  .
\end{gather*}
\end{remark}

\subsection{The modular class in the twisted case}\label{XYZ}
The modular class
of a twisted Poisson structure $(A,\pi,\psi)$ was f\/irst def\/ined in
\cite{yksL} as the $d_{\pi,\psi}$-cohomology class of a section of
$A$ depending on the choice of a volume form (for the case where
$A=TM$, see also \cite{AX}).
Then in \cite{yksW}, it was shown to be the modular class of the morphism
$\pi^\sharp: A^* \to A$, up to a factor 2. Here we
proceed somewhat dif\/ferently.

\begin{definition}\label{modclasstwisted}
The {\emph{modular class}}, $\theta(A ,\pi ,\psi)$,
{\emph{of the twisted Poisson
structure $(\pi,\psi)$ on $A$}} is
\[
\theta(A ,\pi ,\psi) = \frac{1}{2} \, {\rm Mod}   (\pi^\sharp) \ .
\]
A $d_{\pi,\psi}$-cocycle in the modular class
$\theta(A,\pi,\psi)$ is called
a {\emph{modular section}} of $(A,\pi,\psi)$.
\end{definition}

Thus, we take \eqref{onehalf}
to be the def\/inition of the modular class,
and it is clear, because of Theorem \ref{untwisted},
that Def\/inition \ref{modclasstwisted} generalizes Def\/inition~\ref{modclasspoisson}.

\begin{example}\label{example2-sec7} When
$\pi^\sharp : A^* \to A$ is an isomorphism of vector
bundles, the modular class vanishes. For instance, the modular class
of the twisted Poisson structure $(\pi, \psi)$ of Example~\ref{example1-sec7} of
Section~\ref{7.1} above vanishes.
\end{example}

In the case where $A = TM$, the modular class of the morphism
$\pi^\sharp$ is that of the Lie algeb\-roid~$T^*M$. Therefore, Theorem~\ref{theorem1} remains  valid in the twisted case:
The modular class of a~twisted Poisson manifold is
equal to one-half of the modular class of its cotangent bundle Lie
algebroid.

The problem that arises is to determine a section of $A$
which is a representative of the modular class thus def\/ined. Assume
that $A$ is orientable, and let
$\lambda
 \in \Gamma (\wedge^{\rm top} A^{*})$ be a volume form on $A$.
Def\/ine sections $X_\lambda$ and $Y$ of $A$ by
\[
{*_{\lambda} \, X_\lambda = - \partial_{\pi} \lambda}  , \qquad
  {Y = \pi^{\sharp} \iota_{\pi}\psi}    ,
\]
and set $Z_\lambda= X_\lambda +Y$.
In \cite{yksL}, it was proved that $Z_\lambda$ satisf\/ies the relation
\[
\partial_{\pi,\psi, \lambda} - \partial = \iota_{Z_\lambda}   ,
\]
where $\partial_{\pi,\psi, \lambda}
= - *_{\lambda} \circ \, \d _{\pi,\psi} \circ
(*_{\lambda})^{-1}$ and
$\partial = \partial_{\pi} + \underline \partial + \iota_Y$,
the operator $\underline \partial$ being the dual of $\underline
\d$, and also that both operators $\partial_{\pi,\psi,\lambda}$
and $\partial$ are generators of square zero of the
Gerstenhaber bracket $[~,~]_{\pi,\psi}$.
This implies that $Z_\lambda$ is a $\d_{\pi,\psi}$-cocycle.
It is then proved that $Z_\lambda$ is a characteristic cocycle of the
representation $D^{(\partial)}$ of $A^*$ on $\w^{\rm{top}}(A^*)$
 def\/ined by $D^{(\partial)}_\alpha (\mu) = - \alpha \w \partial
 \mu$, for $\alpha \in \Gamma(A^*)$ and $\mu \in
 \Gamma(\w^{\rm{top}}A^*)$.
It was further remarked in \cite{yksW} that, for $\Phi = \pi^\sharp$,
the representation $D^\Phi$ def\/ined by \eqref{repmorphism}
is the square of the representation
$D^{(\partial)}$, $D^{\pi^\sharp} (\mu \otimes \mu) =
D^{(\partial)}
\mu \otimes \mu + \mu \otimes D^{(\partial)} \mu$.
Applying Proposition \ref{charmorphism} to the morphism $\pi^\sharp$
proves the following theorem.
\begin{theorem}\label{X+Y}
The section $Z_\lambda = X_\lambda+Y$
of $A$ is a modular section of $(A, \pi ,\psi)$.
\end{theorem}
We see that in the Poisson case, the modular section
$Z_\lambda$  reduces
to the modular section $X_\lambda$.

See \cite{yksL} for the description of a
unimodular twisted Poisson
structure
on a dense open set of a semi-simple Lie group.
In this example the $3$-form $\psi$ is
the Cartan $3$-form.

\subsection{The regular twisted case}
The formula for the modular class $\theta(A,\pi)$ of a regular Poisson
structure given in Theorem \ref{theoremregular}
is valid for the modular class of a regular twisted Poisson structure
\cite{yksY}:
\[
\theta(A , \pi ,\psi) = (\pi^{\sharp}_{B})^{*} (\theta_{B}) .
\]

\subsection{Non-degenerate twisted Poisson
  structures}\label{7.4}
It is well-known that the inverse of a non-degenerate Poisson
structure, $\pi$, on a Lie algebroid is a symplectic structure,
i.e., a non-degenerate $2$-form, $\omega \in \Gamma(\w^2 A^*)$.
Recall that $\omega$ is def\/ined by $\iota_X(\omega) = - \omega^\flat (X)$,
for all $X \in \Gamma A$, where $\omega^\flat = (\pi^\sharp)^{-1}$.

When a twisted Poisson structure $(A,\pi,\psi)$ is such that $\pi^\sharp$
is invertible, the inverse of $\pi$ is a~$2$-form which is
not $\d_A$-closed, but satisf\/ies
\begin{equation}\label{linear}
\d_A \omega = -\psi  ,
\end{equation}
and, conversely, equation~\eqref{linear} implies that $(\omega^{-1}, \psi)$
is a twisted Poisson structure.
The non-degenerate $2$-form $\omega$ is then
called {\it{twisted symplectic}}.
(See~\cite{SW,yksY},
and, for a generalization of this correspondence, see~\cite{yks2007}.)
It is easy to construct non-degenerate
twisted Poisson structures using this result.

\begin{example}\label{example3-sec7}
We shall show that  Example~\ref{example1-sec7}
of Section~\ref{7.1} arises in this way, when $A$
is the tangent bundle of the manifold, $M$, and $\d_A$ is the de
Rham dif\/ferential, $\d$.
We keep the notations of~\ref{7.1}.
Consider the non-degenerate $2$-form on $M$,
$\omega = (x_3)^{-1} \d x_1 \wedge \d x_2
+ (x_1)^{-1} \d x_3 \wedge \d x_4$. Then $ - \d \omega =
 ((x_1)^{-2} \d x_4 - (x_3)^{-2} \d x_2)
\wedge \d x_1
\wedge \d x_3$.
The inverse of $\omega$ is the non-degenerate bivector $\pi
 = x_3 \, \partial_1 \wedge \partial_2 + x_1
\partial_3 \wedge  \partial_4$. So it follows without computations
that the pair $(\pi, \psi)$, where $\psi
= - \d \omega$, satisf\/ies
\eqref{twisted}.
\end{example}

Let us consider a Lie algebroid, $A$, with a
regular twisted Poisson structure, $(\pi,\psi)$, and let
$B = \pi^{\sharp} (A^{*})$ be the image of $\pi^\sharp$.
The morphism $\pi^\sharp$ def\/ines a skew-symmetric
isomorphism from
$B^{*} = A^{*}/{\rm {Ker}} \, \pi^{\sharp}$ to $B$. We denote the section of $\w^2 B^*$ def\/ined by the
inverse of this isomorphism by
$\omega_{(B)}$. If
$\psi_{(B)}$ is the pull-back
of $\psi$ to the sub-bundle $B$ of $A$, then $\d_B(\omega_{(B)}) = - \psi_{(B)}$.
These considerations lead to the following proposition from
\cite{yksY} which describes a method for constructing
twisted Poisson structures that uses the linearized form
\eqref{linear}
of the twisted Poisson condition \eqref{twisted}.
\begin{proposition}\label{linearize}
 Let $\omega$ be a $2$-form on a Lie algebroid
$A$ whose restriction $\omega|_{B}$ to a Lie subalgebroid $B$ of $A$
is non-degenerate. Then the inverse of $\omega|_B$ and
$\psi = - \d_A \omega$ define a regular twisted Poisson structure on $A$.
\end{proposition}

The bivector
$\pi$ of this twisted Poisson structure is
such that $\pi^\sharp : A^* \to B$ is the composition of
the canonical projection of $A^*$ onto $B^*$ with
the inverse of $\omega|_B$.

\begin{example}\label{example4-sec7}
On $M = \{(x_1, x_2, x_3, x_4, x_5) \in \mathbb R^5 \, | \, x_1 \neq 0 \quad
{\mathrm {and}} \quad  x_5 \neq 0 \}$, consider the $2$-form,
$\omega = x_5 \d x_1 \wedge \d x_2 + x_1 \d x_3 \wedge \d x_4$.
The restriction of $\omega$ to the sub-bundle $B$ of $TM$ generated by
$\partial_1$, $\partial_2$, $\partial_3$ and $\partial_4$
is non-degenerate, with inverse
$\pi = (x_5)^{-1} \partial_1 \wedge \partial_2 +
(x_1)^{-1} \partial_3 \wedge \partial_4$.
Then $(\pi, -\d \omega)$ is a twisted Poisson structure on $M$. Here
$\d \omega =
\d x_1 \wedge \d x_2 \wedge \d x_5  + \d x_1 \wedge
\d x_3 \wedge \d x_4$ and a computation yields
$\frac{1}{2} [\pi, \pi] = - (\wedge ^3 \pi^\sharp) (\d
\omega) = (x_1)^{-2} (x_5)^{-1} \partial_2 \wedge \partial_3 \wedge
\partial_4$.
\end{example}

\subsection[Twisted triangular $r$-matrices]{Twisted triangular $\boldsymbol{r}$-matrices}\label{7.5}

When $A= \g$ is a Lie algebra, and $\pi = r \in \w^2 \g$ and $\psi \in
\w^3 (\g^*)$ satisfy relation \eqref{twisted}, $(r,\psi)$ is called a
{\it{twisted triangular structure}}, and $r$ is called a {\it{twisted
triangular $r$-matrix}}.
The formulas in Proposition \ref{bialgebra} remain valid for the modular
class of a Lie algebra with a twisted triangular $r$-matrix \cite{yksY}:
\[
\theta (\mathfrak{g}, r, \psi ) = -  r^{\sharp}_{({\mathfrak p})}
  (\chi_{{\mathfrak p} , {\rm{Ker}} \, r^{\sharp}}) =
  r^{\sharp}_{({\mathfrak p})} (\chi_{{\mathfrak p},
    \mathfrak{g}/{\mathfrak p}}) \ .
\]

\begin{example}\label{example5-sec7}
Examples of twisted triangular $r$-matrices are given in
\cite{yksY}. They can be constructed as indicated in Proposition
\ref{linearize}.
We quote the simplest example with a non-vanishing modular class.
Let $\p$ be the Lie subalgebra of $\g = {\mathfrak{gl}}_3(\mathbb R)$
with zeros on the third line, which is the Lie algebra of the
group of af\/f\/ine transformations of ${\mathbb R}^2$.
Consider the $2$-cochain on the Lie algebra~$\g$,
$\omega =  e^*_{12} \w e^*_{21} +
e^*_{11} \w
e^*_{13} + e^*_{22} \w e^*_{23}$. The restriction of $\omega$ to $\p$ is
invertible and its inverse is $r = e_{12} \w e_{21} + e_{11} \w e_{13}
+ e_{22} \w e_{23}$. Then $(r, \psi)$, with $\psi = - \d_\g \omega$, is a twisted
triangular structure on $\g$. In this case, $\chi_{\p, \g/\p} = -
(e^*_{11} + e^*_{22})$ and $\theta(\g,r,\psi) = - (e_{13} + e_{23})$.
\end{example}

\subsection{Twisted quasi-Frobenius Lie algebras}
There is a twisted version of the quasi-Frobenius Lie algebras.
A pair $(\omega,\psi)$ is called a {\it {twisted quasi-Frobenius
structure}} on a Lie algebra $\p$ if
$\omega$ is a non-degenerate $2$-cochain on $\p$
and $\psi$ is a~$3$-cocycle on $\p$ such
that $\d_\p\omega = - \psi$.
Then \cite{yksY},
\begin{proposition}
The pair $(\omega,\psi)$ is
a twisted quasi-Frobenius structure on a Lie algebra,
$\g$, if and only if $(\omega^{-1},
\psi)$ is a twisted triangular structure on $\g$.
\end{proposition}
In this correspondence, quasi-Frobenius structures correspond to
triangular structures, so this proposition extends
Proposition \ref{stolin}.

\section[Modular classes of Poisson-Nijenhuis
  structures]{Modular classes of Poisson--Nijenhuis
  structures}\label{section8}

The modular classes of Poisson--Nijenhuis manifolds  were
introduced in
\cite{DF} and further studied in \cite{yksM2006}.
The extension of the def\/inition and properties of the modular classes to
the case of Poisson--Nijenhuis structures on Lie algebroids was carried
out by Caseiro \cite{Caseiro}.
Applications to the
classical integrable systems, the Toda
lattice among them, are to be found in \cite{DF,AD,Caseiro}.

\subsection{Nijenhuis structures}
The \emph{Nijenhuis torsion} of a
$(1,1)$-tensor, $N \in \Gamma(TM \otimes T^*M)$,
on a manifold $M$,
is the $(1,2)$-tensor $[N,N]_{\mathrm{Fr-Nij}}$
def\/ined by
\[
\frac{1}{2} [N,N]_{{\mathrm{Fr}}-{\mathrm{Nij}}}(X,Y)= [NX,NY] - N[NX,Y] - N[X,NY] + N^2[X,Y] ,
\]
for $X$ and $Y \in \Gamma(TM)$.
A $(1,1)$-tensor
is called a \emph{Nijenhuis tensor}
if its Nijenhuis torsion vanishes.
It is well known that powers (with respect to composition) of
Nijenhuis tensors, considered as endomorphisms of the tangent bundle,
are Nijenhuis tensors.
When $N$ is a $(1,1)$-tensor, set
\begin{equation}\label{nij}
[X,Y]_N = [NX,Y] + [X,NY] - N[X,Y] ,
\end{equation}
for $X$ and $Y \in \Gamma(TM)$.

\begin{proposition}
If $N$ is a Nijenhuis tensor on a manifold, $M$,
the bracket $[~,~]_N$ defined by~\eqref{nij} is a Lie bracket on $\Gamma(TM)$,
and the vector bundle $TM$ with anchor $N:TM \to TM$ and Lie bracket
of sections
$[~,~]_N$ is a Lie algebroid.
\end{proposition}

Let us denote the Lie algebroid
$(TM,N,[~,~]_N)$ by $TM_N$.
The Lie
algebroid dif\/ferential of~$TM_N$ on $\Gamma(\w ^\bullet T^*M)$
is $\d_N = [\d,\iota_N]$, where $\d$ is the de Rham dif\/ferential
and $\iota_N$ is the interior product of $1$-forms by the $(1,1)$-tensor
$N$~\cite{yksM}.
The following proposition is from \cite{DF} (see also~\cite{yksM2006}).

\begin{proposition}\label{prop10}
The modular class of the Lie algebroid $TM_N$ is the class of the
$1$-form $\d \Tr N$ in the $\d_N$-cohomology.
\end{proposition}

It follows either from the relation
$[N,N]_{\mathrm{Fr-Nij}} = 0$, or from the fact that
$N$ is the anchor of~$TM_N$, that $N$ is a Lie algebroid morphism from
$TM_N$ to $TM$.
Since the modular class of~$TM$
vanishes,
the modular class of the Lie algebroid $TM_N$ is equal to the modular class
of the Lie algebroid morphism $N$.

\subsection[Poisson-Nijenhuis structures]{Poisson--Nijenhuis structures}

If $\pi$ is a Poisson bivector on $M$,
we denote the Lie algebroid $(T^*M, \pi^\sharp, [~,~]_\pi)$ by
$T^*M_{[0]}$.
A~Nijenhuis tensor $N$ and a Poisson tensor $\pi$ on $M$
are called \emph{compatible}
if $(TM_N, T^*M_{[0]})$ is a~Lie bialgebroid. The pair $(\pi,N)$ is then
called a \emph{Poisson--Nijenhuis structure} (or {\emph{PN-structure}} for
short) on $M$, and $M$ is called a {\emph{Poisson--Nijenhuis manifold}}.

This characterization of Poisson--Nijenhuis structures
\cite{yksLMP} is equivalent to the original de\-f\/i\-nition
\cite{MM,yksM}: $\pi^\sharp N^* = N \pi^\sharp$ and
$C(\pi,N) = 0$, where $N^*: T^*M \to T^*M$
is the dual of $N : TM \to TM$,
and $C(\pi,N)$ is a $(2,1)$-tensor def\/ined by
$C(\pi,N)(\alpha, \beta) = [\alpha,\beta]_{N\pi} -
([\alpha,\beta]_\pi)_{N^*}$, for $\alpha$ and $\beta \in \Gamma(T^*M)$.

The main consequence of the
compatibility of $\pi$ and $N$ is that all the bivectors
in the sequence $\pi_k = N^k \pi$, $k \in \mathbb N$, are Poisson bivectors.
In addition, $(\pi_k)$ constitutes a hierarchy of
pairwise compatible Poisson structures
in the sense that, for any $k$, $k' \in \mathbb N$, $\pi_k + \pi_{k'}$
is a Poisson bivector.

Def\/ine
\[
I_k = \frac{1}{k} \Tr N^k,
\]
for each $k \geq 1$. The
functions $I_k$ are pairwise in involution with respect to any of the
Poisson structures $\pi_\ell$ \cite{MM}, and they are called the \emph{fundamental
  functions in involution} of the Poisson--Nijenhuis structure \cite{MM}.
The corresponding Hamiltonian vector f\/ields,
\[
H^\pi_{I_k} = \frac{1}{k} \pi^\sharp (\d \Tr N^k) ,
\]
are bi-Hamiltonian for $k \geq 2$ with respect to any pair of Poisson structures
in the hierarchy, and they commute pairwise. The sequence of vector
f\/ields $(H^\pi_{I_k})$
is called the {\it{bi-Hamiltonian hierarchy}} of the Poisson--Nijenhuis
structure.

When $(\pi,N)$ is a Poisson--Nijenhuis structure on $M$,
the vector bundle $T^*M$ with anchor
$(\pi_k)^\sharp$ and bracket $[~,~]_{\pi_k}$ is a Lie algebroid which we
denote by $T^*M_{[k]}$. We then denote
$N^{*}$, considered as a map from $T^*M_{[k]}$ to $T^*M_{[k-1]}$, by
$N^{*}_{[k]}$, for each $k \geq 1$.

\begin{proposition}
Let $(\pi,N)$ be a {Poisson--Nijenhuis structure} on $M$. Then, for all $k
\geq 1$, $N^{*}_{[k]} :
T^*M_{[k]} \to T^*M_{[k-1]}$ is a Lie algebroid morphism.
\end{proposition}

\subsection[The modular classes of a Poisson-Nijenhuis
  manifold]{The modular classes of a Poisson--Nijenhuis manifold}  \label{modPN}

\begin{definition}
For $k \geq 1$, the
{\emph{$k$-th modular class}} of a Poisson--Nijenhuis manifold $(M,\pi,N)$
is one-half of the modular class of the Lie algebroid morphism $N^{*}_{[k]}$.
\end{definition}
We denote the $k$-th modular class of $(M,\pi,N)$ by $\theta^{(k)}$.
We recall that, corresponding to each Poisson structure $\pi_k$ on $M$
and consequently to each Lie algebroid $T^*M_{[k]}$, there is associated a
dif\/ferential on $\Gamma(\w^\bullet TM)$, denoted by $\d_{\pi_k}$.
The $k$-th modular class is a class in the f\/irst
cohomology space of the complex
$(\Gamma(\w^\bullet TM), \d_{\pi_k})$.

Assume that $M$ is orientable, and let $\lambda$ be a volume form on $M$.
For each $k \in \mathbb N$, let $X^k_\lambda$ be the modular vector f\/ield
of the Poisson manifold $(M, \pi_k)$ which is a
$\d_{\pi_k}$-cocycle. Then \cite{DF},
\begin{theorem}\label{hierarchy}
Let $(M, \pi,N)$ be an orientable {Poisson--Nijenhuis manifold}. Set
\[
X^{(k)} = X^k_\lambda - N(X^{k-1}_\lambda)  ,
\]
for each $k \geq 1$. Then,

  $\bullet$
$X^{(k)}$ is a $\d_{\pi_k}$-cocycle, independent of the choice of
the volume form, $\lambda$,

 $\bullet$
the class of $X^{(k)}$ in the $\d_{\pi_k}$-cohomology of $\Gamma(\w^\bullet
TM)$ is $\theta^{(k)}$,
the $k$-th modular class of $(M,\pi,N)$,

  $\bullet$
$X^{(k+1)} = N(X^{(k)})$,

  $\bullet$
$X^{(k)} = - \frac{1}{2} H^\pi_{I_k}$.
\end{theorem}

The vector f\/ields $X^{(k)}$ can therefore
be called the {\emph{modular vector
fields of the Poisson--Nijenhuis structure}} $(\pi, N)$ on $M$.
The sequence of modular
vector f\/ields $(X^{(k)})_{k\geq2}$ coincides, up to sign and a factor
$2 $, with the bi-Hamiltonian hierarchy of the
Poisson--Nijenhuis structure.
Under the additional
assumption that $N$ is
invertible, each $X^{(k)}$ is bi-Hamiltonian, not only for
$k \geq 2$, but also for $k=1$ with $I_0= {\rm{ln}} |{\rm{det}} N|$
(see \cite{DF}).

\begin{remark} A proof of Theorem \ref{hierarchy}
is also in \cite{yksM2006},
based on the fact that
$C(\pi,N)=0$ implies
$\iota_\pi(\d_N \d f) = - \frac{1}{2} H^\pi_{I_1} (f)$
which fact is a corollary of a theorem of Beltr\'an and Monterde \cite{BM}.
\end{remark}

The vector f\/ield $X^{(1)} = X^{N\pi}_\lambda - N (X^{\pi}_\lambda)$ is called
\emph{the modular vector field} of $(M,\pi,N)$. Its class in the
$\d_{N\pi}$-cohomology is $\theta^{(1)} = \frac{1}{2} \,
{\rm{Mod}}(N^*)$, called \emph{the modular class} of
$(M,\pi, N)$.
Each~$X^{(k)}$, for $k \geq 2$, is
obtained  from $X^{(1)}$ by the relation,
\[
X^{(k)} = N^{k-1} X^{(1)}  .
\]

Let $TM_{[k]}$ be the Lie algebroid $(TM, N^k, [~,~]_{N^k})$, and let
$N_{[k]}$
be $N$ considered as a Lie algebroid morphism from
$TM_{[k]}$ to $TM_{[k-1]}$.
There is a simple relation between the modular classes of the Lie
algebroid morphisms
$N^{*}_{[k]}: T^*M_{[k]} \to T^*M_{[k-1]}$ and $N_{[k]}: TM_{[k]} \to TM_{[k-1]}$.
\begin{proposition}
When $(\pi,N)$ is a Poisson--Nijenhuis structure on $M$,
\[
{\rm{Mod}}(N^{*}_{[k]}) = -  (\pi^\sharp)^*   ({\rm{Mod}}(N_{[k]}))  ,
\]
for all $k \geq 1$.
\end{proposition}

\subsection[Poisson-Nijenhuis structures on Lie algebroids]{Poisson--Nijenhuis structures on Lie algebroids}

It is straightforward to extend the notion of Poisson--Nijenhuis
manifold to that of Poisson--Nijenhuis structure on a Lie algebroid
\cite{yksM,GU},
more generally on Gelfand--Dorfman complexes
\cite{GD,Dorfman,Vaisman1,Vaisman2}, and the results concerning the
modular classes are readily extended \cite{Caseiro}.
For a~Nijenhuis tensor on a Lie algebroid $A$,
the modular class of the morphism $N$ from $A$ with anchor $a_A \circ
N$ and deformed bracket
$[~,~]_{A,N}$ to $(A,a_A, [~,~]_A)$ is the class of the \mbox{$1$-form}~$\d_A \Tr N$ (cf. Proposition \ref{prop10}).
For a Poisson--Nijenhuis structure on $A$, there is a sequence of
{\it {bi-Hamiltonian modular sections}} (cf.\ Theorem \ref{hierarchy}),
and a multi-Hamiltonian system on the base manifold.
In~\cite{Caseiro}, Caseiro def\/ines a
Poisson--Nijenhuis structure on a Lie algebroid, whose hierar\-chy of
Poisson structures projects onto the
known hierarchy of multi-Hamiltonian structures for the Toda
lattice.

\section{The spinor approach to the modular
 class}\label{section9}
The spinor approach to Poisson geometry stems from Hitchin's
introduction of the generalized geometry using the Courant bracket
on the direct sum $TM \oplus T^*M$ \cite{Hitchin}.
It was developed by Gualtieri \cite{Gualtieri}, and, more recently by
Alekseev, Bursztyn and Meinrenken \cite{Meinrenken,ABM}, and many others.
Here we show how the modular class appears in this approach.

\subsection{The modular class and pure spinors}\label{spinors}
The title of this section will be explained in the next section.

We f\/irst make two def\/initions.
If $\pi$ is a bivector on a vector bundle, $E$, we denote by
$e^{-\iota_\pi}$ the
operator
$\Id -\iota_\pi + \frac{1}{2!} \, \iota_\pi \circ \iota_\pi - \dots +
(-1)^k \frac{1}{k!} (\iota_\pi)^k + \cdots$, whose evaluation on a given form
is a f\/inite sum, and, for $\lambda$ a form of top degree on $E$, let us set
\[
\lambda^{(\pi)} = e^{-\iota_\pi} \lambda  .
\]
Next, we def\/ine, for
$\psi$ a $\d_A$-closed $3$-form on a Lie algebroid~$A$, the
operator on forms,
\[
\d_\psi = \d_A + \epsilon_\psi  ,
\]
where $\epsilon_\psi$
denotes the exterior product by $\psi$. This operator is of
square zero since
$\psi$ is $\d_A$-closed.
Such a closed $3$-form is usually considered as a {\it{background
    $3$-form}}
\cite{KStrobl,SW,Gualtieri}, and $\d_\psi$ is called the {\it{$\psi$-twisted
differential}} associated to $(A,\psi)$. We now show the relationship
of these objects with the modular class.

Let $\pi$ be a bivector, $\psi$ a $\d_A$-closed $3$-form, and $\lambda$
a volume form on an
orientable Lie algeb\-roid,~$A$, and let
$\, X_\lambda$ and $Y$ be def\/ined as in Section~\ref{XYZ},
$*_\lambda X_\lambda = - \partial_\pi \lambda $ and
$  Y = \pi^\sharp \iota_\pi \psi$. Then,
\begin{gather}\label{spinorwopsi}
\d_A (\lambda^{(\pi)})= \iota_{-\frac{1}{2} [\pi,\pi]_A + X_\lambda}
(\lambda^{(\pi)})  ,
\end{gather}
and
\[
\epsilon_\psi (\lambda^{(\pi)}) = \iota_{(\w^3\pi^\sharp) \psi + Y}
(\lambda^{(\pi)})  .
\]
Setting
$  \phi_\pi = \frac{1}{2} [\pi,\pi]_A - (\w^3\pi^\sharp) \psi  $
and $  Z_\lambda= X_\lambda + Y  $,
we obtain
\begin{equation}\label{spinoreq}
\d_\psi (\lambda^{(\pi)})= \iota_{- \phi_\pi + Z_\lambda}
(\lambda^{(\pi)})  .
\end{equation}
These formulas are results of Meinrenken
\cite{Meinrenken}, themselves
based on formulas of Evens and Lu~\cite{EL}. (See also~\cite{AX}.)
Since $\phi_\pi = 0$ is the condition for $(\pi,\psi)$ to def\/ine a twisted
Poisson structure on~$A$, the following theorem is a direct consequence of~\eqref{spinoreq}.

\begin{theorem}\label{theoremspin}
Let $(A, \pi, \psi)$ be an orientable Lie algebroid with a twisted Poisson
structure, and let $\lambda$ be a volume
form on $A$.
The modular section associated to $\lambda$
is the unique section $Z_\lambda$ of $A$ such that
\begin{equation}\label{spinorZ}
\d_\psi (\lambda^{(\pi)})=
\iota_{Z_\lambda}(\lambda^{(\pi)}) .
\end{equation}
An orientable Lie algebroid with a twisted Poisson
structure, $(A,\pi, \psi)$, is unimodular if
and only if there exists a volume form, $\lambda$, on $A$ such that
$\lambda^{(\pi)}$
is a $\d_\psi$-closed form.
\end{theorem}

If, in particular, $\psi = 0$, equation \eqref{spinoreq} reduces to
\eqref{spinorwopsi},
and the theorem states that, for~$(A,\pi)$ an orientable Lie algebroid with a
Poisson structure, $\d_A (\lambda^{(\pi)})=
\iota_{X_\lambda}(\lambda^{(\pi)})$, so that the unimodularity corresponds
to the existence of a volume form, $\lambda$, such that $\lambda^{(\pi)}$ is
$\d_A$-closed.

\begin{remark} \label{remark1-sec9}
In terms of the big bracket, if we again denote the Lie algebroid
structure of~$A$ by~$\mu$,
relation \eqref{spinorZ} becomes
$\{\mu,  \lambda^{(\pi)} \} + \epsilon_\psi  (\lambda^{(\pi)}) = \{Z_\lambda,
\lambda^{(\pi)} \}$,
and in the Poisson case, it takes the simple form,
$\{\mu,  \lambda^{(\pi)} \} =\{X_\lambda, \lambda^{(\pi)} \}$.
Thus, for a Lie algebroid with a twisted Poisson structure, the
modular section associated to the volume form $\lambda$
is the unique section, $Z$, of $A$ satisfying the equation,
\[
\{\mu,  \lambda^{(\pi)} \} + \epsilon_\psi(\lambda^{(\pi)})
=\{Z, \lambda^{(\pi)} \}  ,
\]
and, for a Poisson structure, this equation takes the simple form,
$\{\mu,  \lambda^{(\pi)} \} =\{X, \lambda^{(\pi)} \}$.
\end{remark}

\subsection{Pure spinors}
When the sections of $\w ^\bullet (E \oplus E^*)$ act on the
sections of $\w ^\bullet(E^*)$ by interior and exterior products, $\Gamma
(\w^\bullet E^*)$ appears as a spinor bundle of the Clif\/ford algebra
of $E \oplus E^*$, because of the relation $\iota_X \epsilon_\xi  +
\epsilon_\xi  \iota_X = \langle \xi, X \rangle$, for $X \in \Gamma E$
and $\xi \in \Gamma (E^*)$.
If $\pi$ is a bivector on $E$, its graph
is a Lagrangian sub-bundle of $E\oplus E^*$, i.e., a
maximal isotropic sub-bundle
with respect to the bilinear form def\/ined by the
pairing of $E$ with $E^*$.
A Lagrangian sub-bundle $L$ of $E\oplus E^*$ is determined by a~form,
$\kappa \in \Gamma (\w^\bullet E^*)$,
unique up to a multiplicative constant,
such that $L$ is the annihilator of~$\kappa$
under the
Clif\/ford action  \cite{Chevalley}. Such a form is called a {\it {pure spinor}}.
It can de deduced from a~theorem of
Chevalley \cite{Chevalley} that the graph of~$\pi$ is the annihilator
of the pure spinor~$\lambda^{(\pi)}$. Thus, in
the spinor approach to Poisson and twisted Poisson structures,
the inhomogeneous form $\lambda^{(\pi)}$ is viewed as a pure spinor
def\/ining the graph of $\pi$.

\subsection{Courant algebroids and Dirac structures}

When $A$ is a Lie algebroid, the vector bundle $A \oplus A^*$ becomes
a Courant algebroid with the non-skewsymmetric bracket def\/ined as a
derived bracket of the big bracket by
\[
[ X+\xi, Y +\eta] = \{\{X+\xi, \d_A\},Y+\eta\}  ,
\]
for $X$, $Y \in \Gamma A$ and $\xi$, $\eta \in \Gamma (A^*)$,
and, more generally, if $\psi$ is a $d_A$-closed $3$-form on $A$,
the vector bundle
$A \oplus A^*$ is a Courant algebroid with the non-skewsymmetric bracket,
\begin{equation}\label{twcourant}
[ X+\xi, Y +\eta]_\psi = \{\{X+\xi, \d_\psi\}, Y+\eta\}  ,
\end{equation}
called the {\it{Courant bracket with background}} $\psi$, or the
{\it{$\psi$-twisted Courant bracket}} (see
\cite{Roytenberg,yksderived,yks2005}
and, for further developments, \cite{Terashima,yks2007}).

A Lagrangian sub-bundle of a Courant algebroid whose sections are
closed under the Courant bracket is called a {\it{Dirac structure}}.

By \eqref{spinoreq}, requiring that $\phi_\pi = 0$, the condition for $(\pi,
\psi)$ to be a twisted Poisson structure, is equivalent
to the condition that there exist a section $U$ of $A$ such that
the pure spinor $\lambda^{(\pi)}$ satisfy
\begin{equation}\label{spinor}
\d_\psi( \lambda^{(\pi)})
= \iota_U \lambda^{(\pi)} \ .
\end{equation}
This condition, in turn, is equivalent to
the condition that the annihilator of $\lambda^{(\pi)}$, a
Lagrangian sub-bundle of
$A \oplus A^*$, be closed
under the Courant bracket with background $\psi$ \cite{Gualtieri,Meinrenken}.
Since the annihilator of $\lambda^{(\pi)}$ is the graph of $\pi$,
these considerations imply the following
proposition.
\begin{proposition}
The pair $(\pi,\psi)$ is a twisted Poisson structure on $A$ if and
only if the graph of $\pi$ is a Dirac sub-bundle of $A \oplus A^*$
with the $\psi$-twisted Courant bracket.
\end{proposition}
This proposition is well-known and can be proved without any appeal to
spinors (see \cite{SW, Roytenberg, Terashima}). Here, we emphasize that the
modular f\/ield appears as the section $U$ such that the
$\psi$-twisted dif\/ferential of $\lambda ^{(\pi)}$ satisf\/ies \eqref{spinor},
which ensures that the graph of $\pi$ is closed under the
$\psi$-twisted Courant bracket \eqref{twcourant}.

\pdfbookmark[1]{Appendix}{appendix}
\section*{Appendix: additional references and conclusion}

In this Appendix, we review some recent publications and the direction
of work in progress.

Many lines of research that touch upon the subject of this
survey have been or are being developed.
In the early 90's, the modular vector f\/ields were a tool in the
classif\/ication of Poisson structures, in particular of the quadratic
Poisson structures on the plane \cite{Dufour, LX1992, GMP}.
Work on the classif\/ication of Poisson structures on surfaces using the
periods of a modular vector f\/ield along the singular set of the  Poisson
bivector has been carried out by Radko \cite{Radko1,Radko2}.
Starting in 2000, results have been published on the relation
of the modular
class of a Poisson manifold and the holonomy
of its symplectic foliation
\cite{Ginzburg}, the Godbillon--Vey class of the symplectic
foliation in the regular case \cite{Mikami} and its Reeb class
\cite{AB}, and on the invariance of the modular class under Morita
equivalence \cite{Ginzburg2, BR}. In \cite{Ibort}, Ibort and
Mart\'inez Torres study the unimodularity of
some Poisson manifolds constructed by surgery. In \cite{Mitric},
Mitric and Vaisman compute the modular class of the complete lift of
the Poisson bivector of a manifold to its tangent bundle.
In \cite{Xu2}, Ping Xu determines how the modular class of a Dirac
submanifold of a Poisson manifold is related to that of the ambient
Poisson manifold.
J.-H. Lu recently described the modular classes of
Poisson homogeneous spaces \cite{Lu}, generalizing a result of \cite{ELW}.

At the same time many generalizations of the notion of modular class
have appeared.
The modular classes of Leibniz algebroids
\cite{ILMP, Wade}, of Nambu--Poisson
structures \cite{ILMPNambu, ILLMP},
of Jacobi manifolds \cite{Vaisman3} and
Jacobi algebroids (also
called generalized Lie algebroids) \cite{ILMPgen}
have been def\/ined, their unimodularity implying a
homology/cohomology duality, as well as those of symp\-lectic supermanifolds
\cite{Monterde}.
The modular class of quasi-Poisson $G$-manifolds was def\/ined
by Alekseev et al.\
in \cite{AKM}.

Many further developments are in progress. To name just a few,
work on the various genera\-li\-zations continues, for example
on Jacobi--Nijenhuis algebroids \cite{CC},
while the implications of the properties of the modular classes of
Poisson--Nijenhuis structures on manifolds and on Lie algebroids for
the theory of
integrable systems are being studied
\cite{AD}.
Closely related works are those
of Launois and Richard on Poisson algebras \cite{Launois},
that of Dolgushev \cite{Dolgushev} on
the exponentiation of the modular class into an automorphism of an
associative, non-commutative algebra, quantizing a Poisson algebra,
and that of Neumaier and Waldmann on
the relationship between unimodularity and
the existence of a trace in a quantized
algebra \cite{NW}. The relationship between unimodularity and the
existence of a trace on non-commutative algebras
was already stressed by Weinstein in \cite{Weinstein}
and should be the subject of further work.

There are at least
three directions in which the theory will surely expand: the def\/inition
and study of the modular
classes of Lie groupoids and their behaviour
under Morita equivalence, more generally introducing
stacks and $2$-categories
(following, e.g., \cite{BW, LTX, TZ}),
applications to sigma-models, and the determination of the
relationship with generalized
complex structures. The number of papers on
sigma-models in which
Lie algebroids and twisted Poisson structures appear
is extremely large. In \cite{Bonechi}, Bonechi and Zabzine introduce
the modular class
of the target manifold of the Poisson sigma-model on the sphere:
the unimodularity
appears as a~condition for the
quantization to be well-def\/ined.
There is new work
\cite{LSX} on holomorphic Lie algebroids, and problems such as
the relation between unimodularity  and the
generalized Calabi--Yau mani\-folds~\cite{Hitchin} remain to be studied.

\subsection*{Acknowledgements}
I express my sincere gratitude to the organizers
of the Midwest Geometry Conference in memory of Thomas Branson  which
was held at the
University of Iowa in May 2007.
Thomas Branson and I collaborated on an article in the early 80's, and the
conference was a much appreciated opportunity to honor his memory and
to learn about the considerable development of the
topic in which we were both interested at the time, conformal invariance.

The present article incorporates many results due to
Camille Laurent-Gengoux, Franco Magri, Alan Weinstein and Milen
Yakimov,  which have appeared in our joint publications and
preprints.
I would like to extend
my sincere thanks to them for all they taught me, and I thank Camille,
in particular, for his remarks on an early version of this paper.
My thanks also to Tudor Ratiu for reminding me of the early history
of the bracket of dif\/ferential forms.

This survey takes into account most of the suggestions of the
referees to whom I am most grateful for their
very substantial and constructive comments.

\pdfbookmark[1]{References}{ref}
\LastPageEnding


\begin{thebibliography}{99}

\footnotesize\itemsep=0pt


\bibitem{AB}
Abouqateb A., Boucetta M.,
The modular class of a regular Poisson manifold and the Reeb class of
its symplectic foliation,  \emph{C. R. Math. Acad. Sci. Paris}  {\bf 337}
(2003), 61--66, \href{http://arxiv.org/abs/math.DG/0211405}{math.DG/0211405}.

\bibitem{AM}
Abraham R., Marsden J.E., Foundations of mechanics,
W.A. Benjamin, Inc., New York~-- Amsterdam, 1967.

\bibitem{AM2}
Abraham R., Marsden J.E., Foundations of mechanics,
second edition, revised and enlarged, with the assistance of Tudor
Ra\c tiu and Richard Cushman, Benjamin/Cummings Publishing Co., Inc.,
Reading, Mass., 1978.

\bibitem{AD}
Agrotis M., Damianou P.,
The modular hierarchy of the Toda lattice, \emph{Differential Geom.
Appl.}, to appear,
\href{http://arxiv.org/abs/math.DG/0211405}{math.DG/0701057}.

\bibitem{ABM}
Alekseev A., Bursztyn H., Meinrenken E.,
Pure spinors on Lie groups, \href{http://arxiv.org/abs/0709.1452}{arXiv:0709.1452}.

\bibitem{AKM}
Alekseev A., Kosmann-Schwarzbach Y., Meinrenken E.,
Quasi-Poisson manifolds, \emph{Canad. J. Math.} {\bf 54} (2002), 3--29, \href{http://arxiv.org/abs/math.DG/0006168}{math.DG/0006168}.

\bibitem{AX}
Alekseev A., Xu P.,
Derived brackets and Courant algebroids, unf\/inished manuscript, 2000.

\bibitem{BM}
Beltr\'an J.V., Monterde J.,
Poisson--Nijenhuis structures and the
Vinogradov bracket, {\it Ann. Global Anal. Geom.}  {\bf 12}  (1994), 65--78.

\bibitem{BV}
Bhaskara K.H., Viswanath K., Calculus on Poisson manifolds,
\emph{Bull. London Math. Soc.} {\bf 20}
(1988), 68--72.

\bibitem{Bojowald}
Bojowald M., Kotov A., Strobl T.,
Lie algebroid morphisms, Poisson sigma models, and of\/f-shell closed
gauge symmetries,
\emph{J. Geom. Phys.} {\bf  54}  (2005), 400--426, \href{http://arxiv.org/abs/math.DG/0406445}{math.DG/0406445}.

\bibitem{Bonechi}
Bonechi F., Zabzine M.,
Poisson sigma model on the sphere,  \href{http://arxiv.org/abs/0706.3164}{arXiv:0706.3164}.

\bibitem{Brylinski}
Brylinski J.-L.,
A dif\/ferential complex for Poisson manifolds,  \emph{J. Differential Geom.} {\bf
28}  (1988),  93--114.

\bibitem{BZ}
Brylinski J.-L., Zuckerman G.,
The outer derivation of a complex Poisson manifold,
\emph{J. Reine Angew. Math.} {\bf 506} (1999), 181--189, \href{http://arxiv.org/abs/math.DG/9802014}{math.DG/9802014}.

\bibitem{BR}
Bursztyn H., Radko O.,
Gauge equivalence of Dirac structures and symplectic groupoids,
\emph{Ann. Inst. Fourier (Grenoble)} {\bf 53} (2003), 309--337, \href{http://arxiv.org/abs/math.SG/0202099}{math.SG/0202099}.

\bibitem{BW}
Bursztyn H., Weinstein A.,
Picard groups in Poisson geometry, \emph{Mosc. Math. J.} {\bf  4}
(2004), 39--66, \href{http://arxiv.org/abs/math.SG/0304048}{math.SG/0304048}.

\bibitem{CW}
Cannas da Silva A., Weinstein A.,
Geometric models for noncommutative algebras, {\it Berkeley Mathematics
Lecture Notes}, Vol.~10, Amer. Math. Soc., Providence, RI, 1999.


\bibitem{Caseiro}
Caseiro R., Modular classes of Poisson--Nijenhuis Lie algebroids,
\emph{Lett. Math. Phys.} {\bf 80} (2007), 223--238, \href{http://arxiv.org/abs/math.DG/0701476}{math.DG/0701476}.

\bibitem{CC}
Caseiro R., Nunes da Costa J. M.,
Jacobi--Nijenhuis algebroids and their modular classes,
\href{http://arxiv.org/abs/0706.1475}{arXiv:0706.1475}.

\bibitem{CP}
Chari V., Pressley A.,
A guide to quantum groups, Cambridge University Press,
Cambridge, 1995.

\bibitem{CL}
Chen Z., Liu Z.-J.,
On (co-)morphisms of Lie pseudoalgebras and groupoids,
\emph{J. Algebra} {\bf 316} (2007), 1--31, \href{http://arxiv.org/abs/0710.2149}{arXiv:0710.2149}.

\bibitem{Chevalley}
Chevalley C., The algebraic theory of spinors, Columbia
University Press, New York, 1954;
The algebraic theory of spinors and Clif\/ford algebras, Collected
Works, Vol.~2, Springer-Verlag, Berlin, 1997.

\bibitem{Cornalba}
Cornalba L., Schiappa R.,
Nonassociative star
product deformations for D-brane world-volumes in curved backgrounds,
\emph{Comm. Math. Phys.} {\bf  225}  (2002), 33--66, \href{http://arxiv.org/abs/hep-th/0101219}{hep-th/0101219}.


\bibitem{CDW}
Coste A., Dazord P., Weinstein A., Groupo\"{\i}des symplectiques,
{\it Publications du D\'epartement de Math\'e\-ma\-ti\-ques, Universit\'e Claude
Bernard-Lyon I} {\bf 2A} (1987), 1--62.

\bibitem{Crainic}
Crainic M.,
Dif\/ferentiable and algebroid cohomology, van Est isomorphisms, and
characteristic classes,
\emph{Comment. Math. Helv.} {\bf 78} (2003), 681--721, \href{http://arxiv.org/abs/math.DG/0008064}{math.DG/0008064}.

\bibitem{DF}
Damianou P.A., Fernandes R.L.,
Integrable hierarchies and the modular class, \href{http://arxiv.org/abs/math.DG/0607784}{math.DG/0607784}.

\bibitem{Dolgushev}
Dolgushev V.,
The Van den Bergh duality and the modular symmetry of a Poisson
variety,
\href{http://arxiv.org/abs/math.QA/0612288}{math.QA/0612288}.

\bibitem{Dorfman}
Dorfman I.Ya., Dirac structures and integrability of nonlinear evolution
equations, John Wiley and
Sons, Chichester, 1993.

\bibitem{DL}
Douady A., Lazard M.,
Espaces f\/ibr\'es en alg\`ebres de Lie et en groupes, \emph{Invent. Math.}
{\bf 1} (1966),  133--151.

\bibitem{Drinfeld}
Drinfeld V.G., Hamiltonian structures on Lie groups, Lie bialgebras
and the geometric meaning of classical Yang--Baxter equations,
\emph{Dokl. Akad. Nauk SSSR} {\bf  268}  (1983),  no. 2, 285--287
(English transl.: \emph{Sov. Math. Dokl.} {\bf 27} (1983), no.~2, 68--71).

\bibitem{Drinfeld1989}
Drinfeld V.G.,
Quasi-Hopf algebras,
\emph{Algebra i
Analiz} {\bf  1}  (1989), 114--148 (English transl.: \emph{Leningrad
Math.~J.} {\bf  1}  (1990), 1419--1457).

\bibitem{Dufour}
Dufour J.-P., Haraki A.,
Rotationnels et structures de Poisson quadratiques,
\emph{C. R. Acad. Sci. Paris S\'er. I Math.} {\bf 312} (1991), 137--140.


\bibitem{Ehresmann}
Ehresmann C.,
Cat\'egories topologiques et cat\'egories dif\/f\'erentiables, in Colloque
G\'eom. Dif\/f. Globale (1958, Bruxelles), Centre Belge Rech. Math.,
Louvain, 1959, 137--150.

\bibitem{EL}
Evens S., Lu J.-H.,
Poisson harmonic forms, Kostant harmonic forms, and the $S\sp
1$-equivariant cohomology of $K/T$,  \emph{Adv. Math.} {\bf  142}
(1999), 171--220, \href{http://arxiv.org/abs/dg-ga/9711019}{dg-ga/9711019}.

\bibitem{ELW}
Evens S., Lu J.-H., Weinstein A.,
Transverse measures, the modular class and a cohomology pairing for
Lie algebroids, \emph{Quart. J. Math. Ser.~2} {\bf 50} (1999),
417--436, \href{http://arxiv.org/abs/dg-ga/9610008}{dg-ga/9610008}.


\bibitem{Fernandes0}
Fernandes R.,
Connections in Poisson geometry. I. Holonomy and invariants, \emph{J. Differential Geom.} {\bf 54}  (2000),  303--365,
\href{http://arxiv.org/abs/math.DG/0001129}{math.DG/0001129}.

\bibitem{Fernandes1}
Fernandes R., Lie algebroids, holonomy and characteristic classes,
\emph{Adv. Math.} {\bf  170}  (2002), 119--179, \href{http://arxiv.org/abs/math.DG/0007132}{math.DG/0007132}.


\bibitem{Fernandes2}
Fernandes R.,
Invariants of Lie algebroids, \emph{Differential Geom. Appl.} {\bf 19}
(2003), 223--243, \href{http://arxiv.org/abs/math.DG/0202254}{math.DG/0202254}.

\bibitem{FLS}
Flato M., Lichnerowicz A., Sternheimer D.,
Alg\`ebres de Lie attach\'ees a une vari\'et\'e canonique,
\emph{J. Math. Pures Appl. (9)}  {\bf 54}  (1975), 445--480.

\bibitem{Fuchssteiner}
Fuchssteiner B., Lie algebra structure of
degenerate Hamiltonian and bi-Hamiltonian systems,
\emph{Progr. Theo\-ret. Phys.} {\bf  68}  (1982), 1082--1104.

\bibitem{GD}
Gelfand I.M., Dorfman I.Ya., Schouten bracket and Hamiltonian
operators,
\emph{Funktsional. Anal. i Prilozhen.} {\bf  14}  (1980), no. 3, 71--74
(English transl.:
\emph{Funct. Anal. Appl.} {\bf 14} (1980), no. 3, 223--226).

\bibitem{Gerstenhaber}
Gerstenhaber M.,
The cohomology structure of an associative ring,  \emph{Ann. of Math. (2)}
{\bf 78}  (1963), 267--288.

\bibitem{GG}
Gerstenhaber M., Giaquinto A., Boundary solutions of the classical
Yang--Baxter equation, \emph{Lett. Math. Phys.} {\bf  40}  (1997),
337--353, \href{http://arxiv.org/abs/q-alg/9609014}{q-alg/9609014}.

\bibitem{Ginzburg2}
Ginzburg V.L.,
Grothendieck groups of Poisson vector bundles, \emph{J. Symplectic
  Geom.} {\bf 1} (2001), 121--169, \href{http://arxiv.org/abs/math.DG/0009124}{math.DG/0009124}.

\bibitem{Ginzburg}
Ginzburg V.L., Golubev A., Holonomy on Poisson manifolds and
the modular class,  \emph{Israel J. Math.}  {\bf 122}  (2001), 221--242, \href{http://arxiv.org/abs/math.DG/9812153}{math.DG/9812153}.

\bibitem{GMM}
Grabowski J., Marmo G., Michor P.,
Homology and modular classes of Lie algebroids,
\emph{Ann. Inst. Fourier (Grenoble)} {\bf  56}
(2006), 69--83, \href{http://arxiv.org/abs/math.DG/0310072}{math.DG/0310072}.

\bibitem{GMP}
Grabowski J., Marmo G., Perelomov A. M., Poisson structures: towards a classif\/ication,
\emph{Modern Phys. Lett.~A}  {\bf 8}  (1993), 1719--1733.

\bibitem{GU}
Grabowski J., Urba\'nski P.,
Lie algebroids and Poisson--Nijenhuis structures,
\emph{Rep. Math. Phys.} {\bf  40}  (1997), 195--208, \href{http://arxiv.org/abs/dg-ga/9710007}{dg-ga/9710007}.


\bibitem{Gualtieri}
Gualtieri M., Generalized complex geometry, \href{http://arxiv.org/abs/math.DG/0703298}{math.DG/0703298}.

\bibitem{HM}
Higgins P.J., Mackenzie K., Algebraic constructions in the
category of Lie algebroids, \emph{J. Algebra} {\bf  129}  (1990),  194--230.

\bibitem{Hirota}
Hirota Y.,
Morita invariant properties of twisted Poisson manifolds,
\emph{Lett. Math. Phys.} {\bf 81} (2007), 185--195.

\bibitem{Hitchin}
Hitchin N., Generalized Calabi--Yau manifolds, \emph{Q.~J.~Math.}
{\bf 54}
(2003), 281--308, \href{http://arxiv.org/abs/math.DG/0209099}{math.DG/0209099}.

\bibitem{HY}
Hodges T., Yakimov M.,
Triangular Poisson structures on Lie groups and symplectic reduction, in
Noncommutative Geometry and Representation Theory in Mathematical
Physics, Editors J.~Fuchs  et al.,
\emph{Contemp. Math.} {\bf 391} (2005), 123--134, \href{http://arxiv.org/abs/math.SG/0412082}{math.SG/0412082}.

\bibitem{Huebschmann1}
Huebschmann J.,
Poisson cohomology and quantization,
\emph{J. Reine Angew. Math.} {\bf 408} (1990), 57--113.

\bibitem{Huebschmann3}
Huebschmann J.,
Lie--Rinehart algebras, Gerstenhaber algebras, and Batalin--Vilkovisky
algebras,
{\em Ann. Inst. Fourier (Grenoble)} {\bf 48} (1998), 425--440, \href{http://arxiv.org/abs/dg-ga/9704005}{dg-ga/9704005}.

\bibitem{Huebschmann2}
Huebschmann J.,
Duality for Lie--Rinehart algebras and the modular class,  \emph{J. Reine
Angew. Math.} {\bf  510}  (1999), 103--159, \href{http://arxiv.org/abs/dg-ga/9702008}{dg-ga/9702008}.

\bibitem{ILLMP}
Ib\'a\~ nez R., de Le\'on M., L\'opez B., Marrero J.C.,
Padr\'on E., Duality and modular class of a Nambu--Poisson structure,
\emph{J. Phys. A: Math. Gen.}  {\bf 34}  (2001), 3623--3650, \href{http://arxiv.org/abs/math.SG/0004065}{math.SG/0004065}.

\bibitem{ILMPNambu}
Ib\'a\~nez R., de Le\'on M., Marrero J.C., Padr\'on E.,
Leibniz algebroid associated with a Nambu--Poisson structure,
\emph{J. Phys. A: Math. Gen.}  {\bf 32}  (1999),  8129--8144, \href{http://arxiv.org/abs/math-ph/9906027}{math-ph/9906027}.

\bibitem{ILMP}
Ib\'a\~ nez R., Lopez B., Marrero J.C., Padr\'on E.,
Matched pairs of Leibniz algebroids, Nambu--Jacobi structures and
modular class,  \emph{C. R. Acad. Sci. Paris S\'er. I Math.}  {\bf
333}  (2001), 861--866.

\bibitem{Ibort}
Ibort A., Mart\'{\i}nez Torres D.,
A new construction of Poisson manifolds, \emph{J. Symplectic Geom.} {\bf 2}
(2003), 83--107.

\bibitem{ILMPgen}
Iglesias D., Lopez B., Marrero J.C., Padr\'on E.,
Triangular generalized Lie bialgebroids: homology and cohomology
theories, in Lie Algebroids and Related Topics in dif\/ferential Geometry
(2000, Warsaw), \emph{Banach Center Publ.}, Vol.~54, Polish Acad. Sci.,
Warsaw, 2001, 111--133.

\bibitem{Karasev}
Karas\"ev M.V.,
Analogues of objects of the theory of Lie groups for nonlinear Poisson
brackets, \emph{Izv. Akad. Nauk SSSR Ser. Mat.}  {\bf 50}  (1986),
no.~3, 508--538 (English transl.: \emph{Math. USSR Izv.} {\bf  28} (1987), 497--527).

\bibitem{KStrobl}
Klim\v c\'{\i}k C., Strobl T.,
WZW-Poisson manifolds,  \emph{J. Geom. Phys.}  {\bf 43} (2002),
341--344, \href{http://arxiv.org/abs/math.SG/0104189}{math.SG/0104189}.


\bibitem{yks1992}
Kosmann-Schwarzbach Y.,
Jacobian quasi-bialgebras and quasi-Poisson Lie groups, in
Mathematical Aspects
of Classical Field Theory (1991, Seattle, WA), Editors M.J.~Gotay,
J.E.~Marsden and V.~Moncrief,
\emph{Contemp. Math.} {\bf 132}   (1992), 459--489.

\bibitem{yks1995}
Kosmann-Schwarzbach Y.,
Exact Gerstenhaber
algebras and Lie bialgebroids, \emph{Acta Appl. Math.} {\bf 41}
(1995),  153--165.

\bibitem{yksLMP}
Kosmann-Schwarzbach Y.,
The Lie bialgebroid of a Poisson--Nijenhuis manifold,
\emph{Lett. Math. Phys.} {\bf  38}  (1996), 421--428.

\bibitem{yks1996}
Kosmann-Schwarzbach Y.,
From Poisson algebras to Gerstenhaber algebras,
\emph{Ann. Inst. Fourier (Grenoble)} {\bf 46}
  (1996), 1243--1274.

\bibitem{yksLNP}
Kosmann-Schwarzbach Y.,
Lie bialgebras, Poisson Lie groups and dressing transformations,
in Integrability of nonlinear systems (1996, Pondicherry),
{\it Lecture Notes in Phys.}, Vol.~495, Springer, Berlin, 1997, 104--170.


\bibitem{yks2000}
Kosmann-Schwarzbach Y.,
Modular vector f\/ields and Batalin--Vilkovisky algebras, in
Poisson Geometry (1998, Warsaw), Editors J.~Grabowski and P.~Urba\'nski,
\emph{Banach Center Publ.}, Vol.~51,
Polish Acad. Sci., Warsaw, 2000, 109--129.

\bibitem{yksderived}
Kosmann-Schwarzbach Y., Derived brackets,  \emph{Lett. Math. Phys.} {\bf  69}
(2004), 61--87, \href{http://arxiv.org/abs/math.DG/0312524}{math.DG/0312524}.

\bibitem{yks2005}
Kosmann-Schwarzbach Y., Quasi, twisted, and all that $\ldots$ in Poisson
geometry and Lie algebroid theory, in The Breadth of Symplectic and
Poisson Geometry, Editors J.E. Marsden and T. Ratiu,
{\it Progr. Math.}, Vol.~232, Birkh\"auser,
Boston, MA, 2005, 363--389, \href{http://arxiv.org/abs/math.SG/0310359}{math.SG/0310359}.

\bibitem{yks2007}
Kosmann-Schwarzbach Y.,
Poisson and symplectic functions in Lie algebroid theory, \href{http://arxiv.org/abs/0711.2043}{arXiv:0711.2043}.

\bibitem{yksL}
Kosmann-Schwarzbach Y., Laurent-Gengoux C.,
The modular class of a twisted Poisson structure,  \emph{Trav. Math.} {\bf 16}  (2005), 315--339, \href{http://arxiv.org/abs/math.SG/0505663}{math.SG/0505663}.

\bibitem{yksLW}
Kosmann-Schwarzbach Y., Laurent-Gengoux C., Weinstein A.,
Modular classes of Lie algebroid morphisms, \href{http://arxiv.org/abs/0712.3021}{arXiv:0712.3021}.

\bibitem{yksM}
Kosmann-Schwarzbach Y., Magri F.,
Poisson--Nijenhuis structures,  \emph{Ann. Inst. H. Poincar\'e Phys.
Th\'eor.} {\bf  53}
(1990), 35--81.

\bibitem{yksM2006}
Kosmann-Schwarzbach Y., Magri F.,
On the modular classes of Poisson--Nijenhuis manifolds, \href{http://arxiv.org/abs/math.SG/0611202}{math.SG/0611202}.

\bibitem{yksMo}
Kosmann-Schwarzbach Y., Monterde J.,
Divergence operators and odd Poisson brackets,
\emph{Ann. Inst. Fourier (Grenoble)} {\bf52} (2002), 419--456, \href{http://arxiv.org/abs/math.QA/0002209}{math.QA/0002209}.

\bibitem{yksW}
Kosmann-Schwarzbach Y., Weinstein A.,
Relative modular classes of Lie algebroids,
\emph{C. R. Math. Acad. Sci. Paris} {\bf  341}  (2005),  509--514, \href{http://arxiv.org/abs/math.DG/0508515}{math.DG/0508515}.

\bibitem{yksY}
Kosmann-Schwarzbach Y., Yakimov M.,
Modular classes of regular twisted Poisson
structures on Lie algebroids, \emph{Lett. Math. Phys.} {\bf 80}
(2007), 183--197, \href{http://arxiv.org/abs/math.SG/0701209}{math.SG/0701209}.


\bibitem{KS}
Kostant B., Sternberg S.,
Symplectic reduction, BRS cohomology, and inf\/inite-dimensional
Clif\/ford algebras,
\emph{Ann. Physics} {\bf 176} (1987), 49--113.

\bibitem{Koszul}
Koszul J.-L.,
Crochet de Schouten--Nijenhuis et cohomologie, in
The Mathematical Heritage of \'Elie Cartan (1984, Lyon),
\emph{Ast\'erisque} (1985), 257--271.

\bibitem{Kotov}
Kotov A., Strobl T.,
Characteristic classes associated to $Q$-bundles, \href{http://arxiv.org/abs/0711.4106}{arXiv:0711.4106}.

\bibitem{Kubarski}
Kubarski J.,
The Weil algebra and the secondary characteristic homomorphism of
regular Lie algebroids, in Lie Algebroids and Related Topics in
Dif\/ferential Geometry (2000, Warsaw),
\emph{Banach Center Publ.}, Vol.~54, Polish Acad. Sci.,
Warsaw, 2001, 135--173.

\bibitem{Launois}
Launois S., Richard L., Twisted Poincar\'e duality for some
quadratic Poisson algebras,  \emph{Lett. Math. Phys.}  {\bf 79}  (2007),  161--174, \href{http://arxiv.org/abs/math.KT/0609390}{math.KT/0609390}.


\bibitem{LSX}
Laurent-Gengoux C., Sti\'enon M., Xu P.,
Holomorphic Poisson structures and groupoids, \href{http://arxiv.org/abs/0707.4253}{arXiv:0707.4253}.

\bibitem{LTX}
Laurent-Gengoux C., Tu J.-L., Xu P.,
Chern--Weil map for principal bundles over groupoids, \emph{Math. Z.}
{\bf 255} (2007), 451--491, \href{http://arxiv.org/abs/math.DG/0401420}{math.DG/0401420}.

\bibitem{LR}
Lecomte P., Roger C.,
Modules et cohomologies des big\`ebres de Lie, \emph{C. R. Acad. Sci. Paris
S\'er. I Math.} {\bf  310}
(1990),  405--410, Erratum, 893--894.

\bibitem{LM}
Libermann P., Marle C.-M.,
Symplectic geometry and analytical mechanics,
{\it Mathematics and Its
Applications}, Vol.~35, D. Reidel Publishing Co., Dordrecht, 1987.

\bibitem{Lichnerowicz}
Lichnerowicz A.,
Les vari\'et\'es de Poisson et leurs alg\`ebres de Lie associ\'ees,
\emph{J. Differential Geom.} {\bf  12}  (1977), 253--300.

\bibitem{LX1992}
Liu Z.J., Xu P., On quadratic Poisson structures,
\emph{Lett. Math. Phys.} {\bf 26} (1992), 33--42.


\bibitem{LX}
Liu Z.J., Xu P., Exact Lie bialgebroids and Poisson groupoids,
\emph{Geom. Funct. Anal.} {\bf  6}  (1996), 138--145.

\bibitem{Lu}
Lu J.-H., A note on Poisson homogeneous spaces,
\href{http://arxiv.org/abs/0706.1337}{arXiv:0706.1337}.

\bibitem{LS}
Lyakhovich S.L., Sharapov A.A.,
Characteristic classes of gauge systems, \emph{Nuclear Phys. B} {\bf
703} (2004), 419--453.


\bibitem{M1}
Mackenzie K.,
Lie groupoids and Lie algebroids in dif\/ferential geometry, {\it London
Mathematical Society Lecture Note Series}, Vol.~124, Cambridge University
Press, Cambridge, 1987.


\bibitem{M2}
Mackenzie K.,
General theory of Lie groupoids and Lie algebroids, {\it London
Mathematical Society Lecture Note Series}, Vol.~213, Cambridge University
Press, Cambridge, 2005.


\bibitem{MX}
Mackenzie K., Xu P.,
Lie bialgebroids and Poisson
groupoids, \emph{Duke Math. J.} {\bf  73}  (1994), 415--452.

\bibitem{MM}
Magri F., Morosi C.,
A geometrical characterization of integrable Hamiltonian systems through the theory of Poisson--Nijenhuis manifolds,
{\it Quaderno}, Vol. 19, University of Milano, 1984.

\bibitem{Meinrenken}
Meinrenken E.,
Lectures on pure spinors and moment maps, \href{http://arxiv.org/abs/math.DG/0609319}{math.DG/0609319}.

\bibitem{Menichi}
Menichi L.,
Batalin--Vilkovisky algebra structures on Hochschild cohomology,
\href{http://arxiv.org/abs/0711.1946}{arXiv:0711.1946}.

\bibitem{Mikami}
Mikami K., Godbillon--Vey classes of symplectic foliations,
\emph{Pacific J. Math.}  {\bf 194} (2000), 165--174.

\bibitem{Mitric}
Mitric G., Vaisman I.,
Poisson structures on tangent bundles, \emph{Differential Geom. Appl.}
{\bf 18} (2003), 207--228, \href{http://arxiv.org/abs/math.DG/0108130}{math.DG/0108130}.

\bibitem{Moerdijk}
Moerdijk I., Mr\v cun J.,
Introduction to foliations and Lie groupoids, {\it Cambridge Studies in
Advanced Mathe\-matics}, Vol.~91, Cambridge University Press, Cambridge, 2003.


\bibitem{Monterde}
Monterde J., Vallejo J. A.,
A modular class of even symplectic manifolds,
\emph{Teoret. Mat. Fiz.}  {\bf 132}  (2002), 50--59
(English transl.: \emph{Theoret. and Math. Phys.}  {\bf 132}  (2002), 934--941).

\bibitem{NW}
Neumaier N., Waldmann S.,
Deformation quantization of Poisson structures associated
to Lie algebroids,  \href{http://arxiv.org/abs/0708.0516}{arXiv:0708.0516}.

\bibitem{Palais}
Palais R., The cohomology of Lie rings, in
Proc. Sympos. Pure Math., Vol.~3, Amer. Math.
Soc., Providence, RI, 1961, 130--137.

\bibitem{Park}
Park J.-S., Topological open $p$-branes, in
Symplectic Geometry and Mirror Symmetry (2000, Seoul),
World Sci. Publ., River Edge, NJ, 2001, 311--384.

\bibitem{Pestun}
Pestun V.,
Topological strings in generalized complex space,
\emph{Adv. Theor. Math. Phys.} {\bf 11} (2007), 399--450, \href{http://arxiv.org/abs/hep-th/0603145}{hep-th/0603145}.

\bibitem{Polishchuk}
Polishchuk A., Algebraic geometry of Poisson brackets, in
Algebraic Geometry 7,  \emph{J. Math. Sci. (New York)}  {\bf 84}
(1997), 1413--1444.

\bibitem{Pradines}
Pradines J.,
Th\'eorie de Lie pour les groupo\"{\i}des dif\/f\'erentiables. Calcul
dif\/f\'erentiel dans la cat\'egorie des groupo\"{\i}des
inf\/init\'esimaux,
\emph{C. R. Acad. Sci. Paris S\'er. A-B}  {\bf 264}  (1967), A245--A248.

\bibitem{Radko1}
Radko O.,
A classif\/ication
of topologically stable Poisson structures on a compact oriented
surface,
\emph{J.~Symplectic Geom.} {\bf 1} (2002), 523--542, \href{http://arxiv.org/abs/math.SG/0110304}{math.SG/0110304}.


\bibitem{Radko2}
Radko O., Toward a classif\/ication of Poisson structures on surfaces,
in Quantization, Poisson Brackets and Beyond (2001, Manchester),
Editor T.~Voronov,
\emph{Contemp. Math.} {\bf 315} (2002), 81--88.


\bibitem{Roytenberg}
Roytenberg D.,
Quasi-Lie bialgebroids and twisted Poisson manifolds, {\em
Lett. Math. Phys.}
{\bf 61}  (2002),  123--137, \href{http://arxiv.org/abs/math.QA/0112152}{math.QA/0112152}.

\bibitem{Schwarz}
Schwarz A., Geometry of Batalin--Vilkovisky quantization,
\emph{Comm. Math. Phys.} {\bf  155}  (1993),  249--260, \mbox{\href{http://arxiv.org/abs/hep-th/9205088}{hep-th/9205088}}.


\bibitem{SW}
\v Severa P., Weinstein A.,
Poisson geometry with a
3-form background, in Noncommutative Geometry and String Theory,
Editors Y.~Maeda and S.~Watamura,
{\em Progr. Theoret. Phys. Suppl.}  {\bf 144} (2001), 145--154.

\bibitem{Stolin}
Stolin A.,
On rational solutions of Yang--Baxter
equation for $\sl(n)$, \emph{Math. Scand.} {\bf{69}} (1991), 57--80.

\bibitem{Terashima}
Terashima Y., On Poisson functions, \emph{J. Symplectic Geom.},
to appear.

\bibitem{TZ}
Tseng H.-H., Zhu C.,
Integrating Lie algebroids via stacks, \emph{Compos. Math.} {\bf 142}
(2006), 251--270, \href{http://arxiv.org/abs/math.DG/0405003}{math.DG/0405003}.

\bibitem{Vaintrob}
Vaintrob A.Yu.,
Lie algebroids and homological vector f\/ields,
\emph{Usp. Mat. Nauk} {\bf 52} (1997), 161--162
(English transl.: \emph{Russ. Math. Surv.} {\bf 52} (1997), 428--429).


\bibitem{Vaisman0}
Vaisman I., Lectures on the geometry of Poisson manifolds, {\it Progr.
Math.}, Vol.~118, Birkh\"auser, Basel, 1994.

\bibitem{Vaisman1}
Vaisman I., The Poisson--Nijenhuis manifolds revisited,
\emph{Rend. Sem. Mat. Univ. Politec. Torino} {\bf  52}  (1994), 377--394.

\bibitem{Vaisman2}
Vaisman I., A lecture on Poisson--Nijenhuis structures, in Integrable
Systems and Foliations
(1995, Montpellier), Editors C.~Albert, R.~Brouzet and J.-P.
Dufour, {\it Progr. Math.}, Vol.~145, Birkh\"auser,
Boston, MA, 1997, 169--185.

\bibitem{Vaisman3}
Vaisman I., The BV-algebra of a Jacobi manifold, \emph{Ann. Polon. Math.}
{\bf 73}  (2000),  275--290, \href{http://arxiv.org/abs/math.DG/9904112}{math.DG/9904112}.

\bibitem{Voronov}
Voronov T., Graded manifolds and Drinfeld doubles for Lie
bialgebroids, in Quantization, Poisson Brackets and Beyond (2001, Manchester), Editor T.~Voronov,
\emph{Contemp. Math.} {\bf 315} (2002), 131--168, \href{http://arxiv.org/abs/math.DG/0105237}{math.DG/0105237}.

\bibitem{Wade}
Wade A.,
On some properties of Leibniz algebroids, in Inf\/inite Dimensional Lie
Groups in Geometry and Representation Theory (2000, Washington, DC),
World Sci. Publ., River Edge, NJ, 2002, 65--78.

\bibitem{Weinstein}
Weinstein A., The modular automorphism group
of a Poisson manifold, \emph{J. Geom. Phys.} {\bf{23}} (1997), 379--394.

\bibitem{Xu1}
Xu P.,  Gerstenhaber algebras and BV-algebras in
Poisson geometry,  \emph {Comm. Math. Phys.}
 {\bf 200} (1999),  545--560, \href{http://arxiv.org/abs/dg-ga/9703001}{dg-ga/9703001}.

\bibitem{Xu2}
Xu P., Dirac submanifolds and Poisson involutions,
\emph{Ann. Sci. \'Ecole Norm. Sup. (4)} {\bf 36}  (2003),  403--430, \href{http://arxiv.org/abs/math.SG/0110326}{math.SG/0110326}.

\bibitem{Zabzine}
Zabzine M.,
Lectures on generalized complex geometry and supersymmetry,
\emph{Arch. Math. (Brno)} {\bf 42}  (2006),  suppl., 119--146, \href{http://arxiv.org/abs/hep-th/0605148}{hep-th/0605148}.



\end{thebibliography}
\end{document}